\documentclass{article}
\usepackage{IREP1}
\usepackage[pdftex]{graphicx}
\usepackage{epstopdf}
\usepackage{fancyhdr}
\usepackage{amsmath}
\usepackage{amssymb}
\def\no{\noindent}

\title{
\LARGE \bf A  formula for damping interarea  oscillations\\ with generator redispatch} 
\author{Sarai Mendoza--Armenta$^{1,2}$~~~~~~~~~~~~Ian Dobson$^1$\\[6pt]
\leftline{\normalsize \!\!\!\!(1) Electrical and Computer Engineering Dept., Iowa State University,
Ames IA 50011 USA, dobson@iastate.edu~~~~~~~~~}\\
\leftline{\normalsize \!\!\!\!(2) Instituto de F\'{\i}sica y Matem\'aticas, Universidad Michoacana, Morelia, Michoac\'an, M\'exico, sarai.mendoza@gmail.com}\\~\\~}
\begin{document}
\fancyhead[c]{\textnormal{ 2013 IREP Symposium--Bulk Power System Dynamics and Control--IX (IREP), August 25-30, 2013, Rethymon, Greece}}
\renewcommand{\headrulewidth}{ 0.0pt}
\renewcommand{\footskip}{ 10pt}
 \fancyfoot[L]{~\\[-20pt]Preprint \copyright 2013 IEEE. \small Personal use of this material is permitted. Permission from IEEE must be obtained for all other uses, in any current or future media, including reprinting/republishing this material for advertising or promotional purposes, creating new collective works, for resale or redistribution to servers or lists, or reuse of any copyrighted component of this work in other works.}
 \fancyfoot[C]{~ }
\maketitle
\begin{abstract}
\looseness=-1
We derive a  new formula for the sensitivity of electromechanical oscillation damping with respect to generator 
redispatch.  The formula could lead to some combination 
of observations, computations and heuristics to more effectively  damp interarea oscillations.
\end{abstract}

\section{Introduction}
Large-scale power systems have multiple low frequency  and lightly damped electromechanical modes of oscillation. 
If the damping of these modes becomes too small or positive, then the resulting oscillations can cause equipment
damage, malfunction or blackouts.
The oscillations can appear for large or unusual power transfers, and may become more frequent 
as power systems experience  greater variability of loading conditions.
 Practical limits for power system security often require sufficient damping of oscillatory modes  \cite{Cigre,Rogersbook} 
and power transfers on tie lines are sometimes limited by oscillations \cite{Cigre,IEEE90,IEEE94,ChungPS04}.

There are several approaches to suppressing low-frequency oscillations, including 
limiting power transfers, installing closed loop controls, and taking control actions such as redispatching generation.
In this paper we do analysis to underpin the suppression of oscillations via redispatch of generation.
 \renewcommand{\footskip}{20pt}
 \pagestyle{plain}
Changes in generator dispatch change the oscillation damping 
by exploiting the nonlinearity of the power system: changing 
the dispatch changes the operating equilibrium and
hence the linearization of the power system about that 
equilibrium that determines the oscillatory
modes and their damping.
The use of generator dispatch to damp oscillations has been demonstrated and there are 
several previous approaches:
\begin{enumerate}
\item There are heuristics in terms of the mode shapes for the redispatch for some simpler grid structures \cite{FischerIREP,FischerPPT}.
\item There are exact computations of the sensitivity of the damping from a dynamic power grid model \cite{DobsonCDC92,DobsonPSERC99,NamPS00,WangJZU08}.
The formulas from these computations require Hessians and left eigenvectors of the mode shapes, or derivatives of eigenvectors.
\item The effective generator redispatches can be determined by repetitive computation of eigenvalues  of a dynamic power grid model to give numerical sensitivities \cite{ChungPS04,mangoPESGM,mangoReport,DiaoPSCE11}.
\end{enumerate}
The requirement in approaches 2 and 3 of 
a  large scale power system dynamic model poses some difficulties. 
It is challenging to obtain validated models of generator dynamics over a wide area and particularly difficult to determine dynamic load models.
Moreover, in approach 2, it does not seem feasible to estimate the left eigenvectors of the mode shapes or derivatives of eigenvectors from measurements.

The new formula for modal sensitivity we derive largely depends on power system quantities that can, at least  in principle, be observed
from measurements.
In particular, the formula shows that the first order effect of a generator redispatch largely depends on the mode shape and the power flow.
(The assumed equivalent generator dynamics only appear as a factor common to all redispatches.)
The mode shape (right eigenvector of the oscillatory mode of interest) is  to some 
considerable extent obtainable from power system measurements \cite{JonssonPS04,TrudnowskiPS08,ChaudhuriPS11,DosiekPS13}.
As a general goal, we would like to 
move towards approaches that take  more advantage of synchrophasor measurements, and are less dependent on 
detailed wide area power system dynamic models that are hard to obtain.
While we have not yet proved that the formula can be the basis for doing this,
the formula does open up this promising possibility.

Another possible approach would be to use the new formula to gain insights into oscillation damping 
by generator redispatch that can be expressed and applied as heuristics.  In fact our 
analytic work is inspired by the heuristics by Fisher and Erlich in \cite{FischerIREP,FischerPPT}, and 
a general goal is to confirm, refine and extend their heuristics by supplying an analytic basis for the heuristics.

A key barrier to better understanding  and computing or deriving 
heuristics for generator redispatch to damp oscillations 
has been the difficulty of  the analysis. 
In this paper we are able to  combine several new and old methods of analysis to 
derive a new formula for the sensitivity of the oscillation damping and frequency with respect to
generator redispatch. Most of the paper is devoted to deriving the new formula,
but we also give some special cases and simple examples to begin the process of 
understanding the formula and how it might be applied. The paper results will also appear in 
Spanish as part of PhD thesis \cite{MendozaArmentaPhD}.

\section{Notation}

We use the following notation and definitions.
All quantities are in per unit unless otherwise stated.

\noindent
  \begin{tabular}{ @{}ll @{}}
 $n$&number of buses\\
 $m$&number of generator buses;\\
 &buses 1,2,\, ...\,,$m$ are the generator buses\\
 $\ell$&number of transmission lines\\[4pt]
  $\delta_i$&  voltage angle of bus $i$\\
  $\delta$&  vector of  voltage angles $(\delta_{1},\delta_{2},...,\delta_{n})^T$\\
   $V_i$&  voltage magnitude at bus $i$\\
 $V$&vector of voltage magnitudes  $(V_{m+1},V_{m+2},...,V_{n})^T$\\
 &\hspace{3ex} at load buses\\ 
 $z$&state vector $(\delta,V)^T$\\
 $V_i^{\rm ln}$&  ${\rm ln} V_i$;  logarithm  of voltage magnitude at bus $i$\\[4pt]
 
 $P_i$&net real power injection at bus $i$\\ 
  $Q_i$&net reactive power injection at bus $i$\\ [4pt]
  $\omega_0$ &
nominal frequency in rad/s\\
$h_i$&  inertia at bus $i$ in seconds\\
$m_i$& $2h_i/\omega_0$\\
$M$& diagonal matrix ${\rm diag}\{m_1,m_2,...,m_{{2n-m}}\}$\\
$d_i$ &  damping coefficient at bus $i$ in seconds\\
$D$& diagonal matrix ${\rm diag}\{d_1,d_2,...,d_{{2n-m}}\}$\\[4pt]
$b_{ij}$&  imaginary part of $ij$ element of  bus admittance \\
 & matrix (for $i\ne j$, $b_{ij}$ is the absolute value of the\\&susceptance
 of the line joining bus $i$ and bus $j$; \\
 &$b_{ii}$ is the sum of the susceptances incident on bus $i$.)\\[4pt]
$b_{k}$&  absolute value of the susceptance of  line $k$ \\
$p_{k}$&  real power flow in  line $k$  defined in (\ref{p})\\
$q_{k}$&  part of  reactive power flows in line $k$  defined in (\ref{q})\\[4pt]
$R$&scalar potential energy function defined in (\ref{R})\\
  $L$&a weighted Laplacian matrix;  Hessian of $R$ \\
  $L^{\dagger}$&Matrix pseudo-inverse of $L$\\
   $Q$& $M\lambda^2+D\lambda+L$; a quadratic matrix function of $\lambda$\\
  $x$& eigenvector of $Q$ \\
    \end{tabular}
    
The bus-line incidence matrix $A$ is  defined by
\begin{align}
A_{ik}=\left\{
\begin{array}{rl}
1& \text{if bus $i$ is the sending end of line $k$,}\\
-1 & \text{if bus $i$ is the receiving end of line $k$,}\\
 0 & \text{otherwise.}
\end{array} \right.\notag
\end{align}
The matrix of absolute values of the entries of $A$ is written as $|A|$.
That is,
\begin{align}
|A_{ik}|=\left\{
\begin{array}{ll}
1 & \text{if bus $i$ is the sending or receiving end of line $k$,}\\
 0 & \text{otherwise.}
\end{array} \right.\notag
\end{align}

The angle across line $k$ is defined as
\begin{align}
\theta_k=\sum_{r=1}^n \!A_{rk}\delta_r=\left\{
\begin{array}{ll}
\delta_i -\delta_j & \text{if bus $i$ is sending end of line $k$},\\
\delta_j -\delta_i &\text{if bus $i$ is receiving end of line $k$},\\
\end{array} \right.\notag
\end{align}
and we write $\theta=(\theta_1,\theta_2,....,\theta_{\ell})^T$ 
for the vector of angles across all the lines.

The new voltage coordinate   for line $k$ that connects bus $i$ to bus $j$ 
is defined by 
\begin{align}
\nu_k & =\sum_{r=1}^n |A_{rk}|\ln{V_r} = \ln{(V_iV_j)},
\notag
\end{align}
and we write $\nu=(\nu_1,\nu_2,....,\nu_{\ell})^T$ 
for the vector of voltage coordinates for  all the lines.

\noindent
  \begin{tabular}{ @{}ll @{}}
$z'$&state vector  $(\theta,\nu)^T$ (``dashed" coordinates)\\
$h$& $z'=h(z)$; transforms  undashed to dashed coordinates\\
$H$&Jacobian of coordinate change $h$\\
$x'$&eigenvector of $Q$ in dashed coordinates $z'$ \\
 \end{tabular}

\section{Power System Model}

We model the generators with simple swing dynamic equations and also consider the real and reactive  power balance algebraic 
equations of the network.  The transmission lines are lossless.

The real power balance equations for all buses are
\begin{align}
\label{active balance with A}
m_i\ddot{\delta}_i + d_i\dot{\delta}_i +\sum_{j\sim i}b_{ij} V_i V_j&\sin(\delta_i-\delta_j)  =
P_i, \notag\\
&  i=1,2,\ldots,n.
\end{align}
The notation $j\sim i$ means that the summation is over all buses $j$ connected to bus $i$, excluding $i$.

The reactive power balance equations for all load buses  are
\begin{align}
\label{reactive balance with A}
-\sum_{j\sim i}b_{ij} V_j\cos(\delta_i-\delta_j) = &\, \frac{Q_i}{V_i}+b_{ii}V_i,\notag\\& i=m+1,\ldots,n.
\end{align}
Note that the reactive power balance equations (\ref{reactive balance with A}) have been divided by 
the bus $i$ voltage magnitude $V_i$  \cite{DeMarcoWassnerCCA95}.

The model (\ref{active balance with A}) and 
(\ref{reactive balance with A}) is differential-algebraic equations written in terms of state variables
$z=(\delta,V)$ of the bus voltage angles $\delta$ and the voltage magnitudes of the load buses $V$.
(The generator voltage magnitudes $V_1,V_2,...,V_m$ are assumed to be constant.)
The coupling of the machine angle dynamics of (\ref{active balance with A}) into the voltages of 
(\ref{reactive balance with A}) is emphasized in \cite{VanfrettiIREP10}.

There are two types of buses:
\begin{enumerate}

\item {\bf Generators.} Generator $i$ is assumed to have constant voltage magnitude $V_i$ and 
 the overall effect of its dynamics described by the swing equation.

\item {\bf Loads.} 
Load $i$ can be modeled as a constant power load with real power $P_i<0$ and reactive power $Q_i$.
Moreover, $m_i=d_i=0$.
However, if desired, it is straightforward to model the frequency dependence of real power with $d_i\ne0$, 
and to allow the reactive power $Q_i$ to be a function of $V_i$.
Connecting buses are a special case of load buses with $m_i=0$, $d_i=0$, $P_i=0$, and $Q_i=0$.

\end{enumerate}

Now we rewrite (\ref{active balance with A}) and 
(\ref{reactive balance with A}) in
terms of the partial derivatives of the scalar potential energy function
\begin{align}
\label{R}
R & = -\sum_{\substack{i,j\\i\ne j, i\sim j}}b_{ij}V_iV_j\cos(\delta_i - \delta_j)\notag\\
& \qquad - \sum_{i=1}^n(P_i\delta_i + {\textstyle \frac{1}{2}}b_{ii}V_i^2 + Q_i\ln V_i). 
\end{align}
The first summation in (\ref{R}) is over all the lines.
$R$ is well-known from energy function approaches to power systems \cite{BergenHillPAS81,AraposthasisCAS82,NarasimhamurthiCAS84,TsolasCAS85,Pai89,DeMarcoWassnerCCA95}.
Then  (\ref{active balance with A}) and (\ref{reactive balance with A})
can be rewritten as
\begin{align}
\label{power balance}
 m_i\ddot{\delta}_i + d_i\dot{\delta}_i+\frac{\partial R}{\partial \delta_i} & = 0,\quad i=1,2,\ldots,n,\\
\label{reactive balance}
\frac{\partial R}{\partial V_i} &= 0,\quad  i=m+1,\ldots,n.
\end{align}

The model is differential-algebraic equations with parameters $P_i$ representing the generator power injections 
that we seek to change to best damp the oscillations.
Note that the parameters $P_i$ do not appear explicitly  in the Jacobian  of 
(\ref{power balance}-\ref{reactive balance});
the  mechanism of the damping is that changing $P_i$ changes the operating point at which the Jacobian 
is evaluated and, since the power system is nonlinear, changes the eigenvalues and the damping.
The generator redispatch damping  is an open-loop control exploiting nonlinearity.

\section{Linear Stability Analysis}

The dynamics of the system is described by a set of nonlinear differential-algebraic
equations, and we will apply linear stability analysis to compute
the electromechanical nodes of the system. The state vector $z=(\delta,V)$. Linearizing 
(\ref{power balance}-\ref{reactive balance}) 
and evaluating at the operating equilibrium  $z_*$, we have
\begin{align}
\label{linearization real power balance}
m_i\ddot{\Delta \delta_i} + d_i\dot{\Delta \delta_i} + \sum_{j=1}^{2n-m} L_{ij}\Delta z_j  & = 0, \ i=1,2,\ldots,n,\\ 
\label{linearization reactive power balance}
\sum_{j=1}^{2n-m} L_{ij}\Delta z_j  & = 0,\  i=m+1,\ldots,n.
\end{align}
The linearized deviations from the equilibrium are written as $\Delta z=(\Delta\delta,\Delta V)$.
The weighted Laplacian matrix $L$ is defined as the Hessian of $R$ evaluated at the equilibrium:
\begin{align}
\label{laplacian}
L=\frac{\partial^2 R}{\partial z^2}\Big|_{z_*}
=\left.
\begin{pmatrix}\displaystyle
\frac{\partial^2 R}{\partial \delta^2} & \displaystyle\frac{\partial^2 R}{\partial \delta\, \partial V}\\[5pt]
\displaystyle\frac{\partial^2 R}{\partial V \partial \delta} & \displaystyle\frac{\partial^2 R}{\partial V^2}
\end{pmatrix}\right|_{z_*}
\end{align}
$L$ is a symmetric $(2n-m)\times (2n-m)$ matrix. We can use
the matrices $M={\rm diag}\{m_1,m_2,...,m_{2n-m}\}$ and $D={\rm diag}\{d_1,d_2,...,d_{2n-m}\}$ to rewrite 
(\ref{linearization real power balance}) and (\ref{linearization reactive power balance})
in matrix form as 
\begin{align}
\label{xnce}
M\ddot{\Delta z} + D\dot{\Delta z} + L\Delta z=0.
\end{align}

Now, following \cite{MalladaCDC11}, we define the quadratic matrix function
\begin{align}
\label{quadratic}
Q(\lambda) = M\lambda^2 +D\lambda + L.
\end{align}
$Q$ is a symmetric complex matrix.
We consider the quadratic eigenvalue problem
of finding  $(\lambda,x) \in \mathbb{C}\times\mathbb{C}^{2n-m}$
such that
\begin{align}
\label{Qx}
Q(\lambda)x=0.
\end{align}
$x$ is the right eigenvector
associated to the eigenvalue $\lambda$.

We assume throughout the paper that the Jacobian evaluated at the operating equilibrium $z_*$ 
has no zero eigenvalues except for those associated with the uniform increase of all the angles.
Moreover, we assume throughout the paper that the eigenvalue $\lambda$ is nonresonant
(algebraic multiplicity one).  These generic assumptions
ensure that $\lambda$ and $x$ (suitably normalized) are locally smooth functions of system parameters.
Although the calculations apply to any nonresonant eigenvalue,
it is convenient to consider in the sequel a particular complex eigenvalue $\lambda$
corresponding to an interarea oscillation.
Then $x$ is a complex eigenvector.

Since $M$, $D$, $L$, and $Q$ are
symmetric, the left eigenvector is the row vector $x^T$;
i.e., $x^TQ=0$.
 
We have 
\begin{align}
\label{Q}
\sum_{i,j}Q_{ij}x_ix_j=x^TQx=0.
\end{align}
It might seem more natural at this point to write 
${\bar x}^T Qx=0$ instead of (\ref{Q}), but nevertheless 
it is important to proceed with  (\ref{Q}).

Using a state vector $( \omega_g,\delta)$, where 
$\omega_g$ is the vector of generator speeds, the second order parts of (\ref{linearization real power balance}) can be rewritten 
as first order differential-algebraic equations, and an extended Jacobian $J$ obtained.
The appendix 
 shows that the finite generalized eigenvalues of $J$ 
are the same as the eigenvalues of the quadratic eigenvalue problem and that 
the eigenvectors also correspond.
In particular, eigenvector $x$ of the quadratic eigenvalue problem with 
eigenvalue $\lambda$ corresponds exactly 
to a right generalized eigenvector $(\lambda x_g,x)$ of $J$, where $x_g$ is the vector of components of $x$ corresponding to generator angles.
Thus there is a direct relationship between the eigenvector $x$ and eigenvalue $\lambda$ of (\ref{Qx}) and the conventional 
right eigenvector eigenstructure of the differential-algebraic equations.

\section{Eigenvalue Sensitivity}
In this section we compute the sensitivity of
the electromechanical modes of the system starting from  (\ref{Q}).
We  suppose that $M$ and $D$ are 
constant matrices.\newline

Computing the differential 
of (\ref{Q}) and using (\ref{Qx}),
\begin{align}
0&=d(x^TQx) =
(dx^T)Qx +
x^T(dQ)x +
x^TQ\, dx\nonumber \\
& =
x^T\!(dQ)x =
x^T\!\left(d(\lambda^2 M + \lambda D + L)\right)x
\nonumber\\
& =
\label{dQ}
2\lambda x^TMx \,d \lambda+
x^TDx \,d \lambda+ 
x^T(dL)x.
\end{align}
Then solving for $d\lambda$ gives the complex equation
\begin{align}
 \label{dlambda0}
d\lambda =
-\frac{x^TdL\,x}
{2\lambda x^TMx + x^TDx}.
\end{align}
The linearization (\ref{dlambda0}) captures the first order sensitivity
of a mode. 
Several conclusions may be drawn from (\ref {dlambda0}).
The sensitivity of a mode
depends on its associated eigenvector $x$, 
but not on changes of the eigenvector. This  means that to predict to first order the changes of 
a mode one does not need to take into account the 
variation of its mode shape. 
Moreover, the change in a mode is proportional to the changes in the network $dL$ 
caused by the generator redispatch.
It is convenient to define the complex number 
\begin{align}
\label{alpha}
\alpha=2\lambda x^TMx + x^TDx.
\end{align}
Then (\ref{dlambda0}) becomes
\begin{align}
 \label{dlambda}
d\lambda =
-\frac{x^TdL\,x}{\alpha}.
\end{align}
Given a particular mode of interest $\lambda$,
the generator dynamical parameters $M$ and $D$ 
and the eigenvalue $\lambda$ only appear in (\ref {dlambda}) as the complex factor $\alpha$
in the denominator that is the same for all redispatches.




As (\ref{dlambda}) depends on  the differential $dL$
of the Laplacian, the
next subsections will, after introducing new coordinates,  focus in calculating $dL$.


\subsection{New Coordinates Related to the 
Transmission Lines }
\label{new coordinates}

We know that  the change in  the 
eigenvalue is proportional to the change $dL$
 of the Laplacian of the 
system. $L$ carries all the information
of the network; i.e., $L$ describes aspects of 
the power flow through every transmission
line of the network, so this suggests 
computing $L$ with coordinates that are 
related to the transmission lines of the network,
instead of coordinates that are related
to the buses. For this reason we define
new dashed coordinates $z'=(\theta,\nu)$. $z'$ is
a vector of size $2\ell$, where $\ell$ is the 
number of transmission lines in the 
network. For line $k$ joining buses $i$ and $j$,
the variables $\theta_k$ and $\nu_k$ are 
defined by
\begin{align}
\label{thetak}
\theta_k & =\sum_{r=1}^nA_{rk}\delta_r =\delta_i - \delta_j,\\
\label{nuk}
\nu_k & =\sum_{r=1}^n |A_{rk}|\ln{V_r} = \ln{(V_iV_j)}.
\end{align}
Equations (\ref{thetak}) and (\ref{nuk}) are a nonlinear change of coordinates
\begin{align}
\label{h}
(\theta,\nu)&=h(\delta,V),\quad\mbox{or} \\
z'&=h(z). 
\end{align}
The Jacobian of the coordinate change $h$ is a 
${2\ell \times (2n-m)}$ matrix written as
\begin{align}
\label{H}
H=\frac{\partial h}{\partial z},
\end{align}
with the entries 
\begin{align}
\label{Hthetak}
H_{ki} &=\left\{
\begin{array}{ll}
A_{ik}, & i=1,\ldots,n\\
0, & i=m+1,\ldots,n\\
\end{array} \right\},\quad k=1,\ldots,\ell,\\
\label{Hnuk}
H_{ki} &=\left\{
\begin{array}{ll}
0, & i=1,\ldots,n\\
{\displaystyle \frac{|A_{ik}|}{V_i}}, & i=m+1,\ldots,n\\
\end{array} \right\},\quad k=\ell+1,\ldots,2\ell.
\end{align}
Note that (\ref{Hnuk}) depends on the magnitude
of the load voltage $V_i$. The eigenvectors $x=(x_{\delta},x_{V})$
transform to eigenvectors $x'=(x'_{\theta},x'_{\nu})$ according to
\begin{align}
\label{xdash}
x'=Hx.
\end{align}
That is, for $k=1,\ldots,\ell$,
\begin{align}
\label{x dash thetak}
x'_{\theta_k}&=\sum_{r=1}^{n}A_{rk}x_{\delta_r}\notag\\&=\left\{\begin{array}{ll}
x_{_{\delta_i}}-x_{_{\delta_j} }& \text{if bus $i$ is sending end of line $k$},\\
x_{_{\delta_j}}-x_{_{\delta_i} }&\text{if bus $i$ is receiving end of line $k$}.\\
\end{array}\right. \\
\label{x dash nuk}
x'_{\nu_k} &=\sum_{r=m+1}^{n}\frac{|A_{rk}|}{V_r}x_{_{V_r}}\notag\\
&=\left\{\begin{array}{ll}
\displaystyle\frac{x_{\scriptscriptstyle {V_i}}}{V_i} + \frac{x_{_{V_j}}}{V_j} 
&\text{if line $k$ joins load bus $i$ to load bus $j$}, \\[4mm]
\displaystyle\frac{x_{\scriptscriptstyle {V_i}}}{V_i}&
\text{\hspace{-8mm}if line $k$ joins load bus $i$ to generator bus $j$}. 
\end{array}\right.
\end{align}
(In (\ref{x dash nuk}), we neglect the case that two generator buses are joined by a line.
This case can be excluded by combining together models for generators behind a common step-up transformer,
and modeling the transformer as a bus.)

In most cases, the new coordinates are overdetermined or
redundant; that is, the system has more line coordinates
than the number of independent voltage angles and magnitudes.\footnote{
Consider a tree network composed by $n$ buses.
The tree network has $n-1$ lines, so the number of line angle
coordinates is equal to the number of independent bus 
voltage angles. If the network has just one generator, the number 
of line voltage coordinates is equal to the number of
voltage magnitude variables. However, 
if the network has more than one generator, 
the number of line voltage coordinates 
is larger than the number of bus voltage magnitudes.
In meshed networks the new line coordinates are always
overdetermined because the number of lines is larger
that the number of independent voltage angles and magnitudes.
}
This does not affect the derivation of the formula in this paper, which applies generally, but 
the dependencies between the line coordinates should be 
kept in mind in future work applying the formula.\footnote{We note the approach 
in \cite{BergenHillPAS81} of using the line angle coordinates corresponding 
to a spanning tree of lines.}

The Laplacian matrix in the new coordinates is
\begin{align}
\label{Ldash}
L'_{ij}=
\frac{\partial^2 R}{\partial z'_i \partial z'_j},
\end{align}
$L'$ is a $2\ell \times 2\ell$ matrix.
The partial derivatives transform according to
\begin{align}
\label{derivatives transformation}
\frac{\partial}{\partial z_i}=\sum_{k=1}^{2\ell}H_{ik}\frac{\partial}{\partial z'_k}.
\end{align}
%
Then $L$ and $L'$ are related by 
\begin{align}
\label{L and Ldash}
L_{ij}&=
\frac{\partial^2 R}{\partial z_i \partial z_j}\notag\\
&=\sum_{h,k=1}^{2\ell}
H_{ih}H_{kj}\frac{\partial^2 R}{\partial z'_h \partial z'_k}
=\sum_{h,k=1}^{2\ell}
H_{ih}H_{kj}L'_{hk},
\end{align}
or 
\begin{align}
\label{L and Ldash2}
L&=H^TL'H.
\end{align}
Then
\begin{align}
\label{dL1}
dL& =d(H^TL'H)= H^T(dL')H + 2H^TL' dH.
\end{align}

\subsection{Computing dL}
In this section the goal is to compute the right hand side of
expression (\ref{dL1}). To do it we have to express $R$ from 
(\ref{R}) in terms of $z'$, but note that $R$ is naturally
composed by one part related to the transmission lines
 and another part that refers to
the buses of the system.  As $z'$ is related to
the transmission lines, this suggests to express just the first
part of $R$ in terms of $z'$ and to keep the second part in terms
of the bus coordinates; i.e.,
\begin{align}
\label{R def}
R & = R_{\rm line}+R_{\rm bus},
\end{align}
where
\begin{align}
\label{R lines coordenates}
 R_{\rm line}&=-\sum_{\substack{i,j\\i\ne j, i\sim j}}b_{ij}V_iV_j\cos(\delta_i - \delta_j),\\
\label{R buses coordenates}
R_{\rm bus}&=-\sum_i^n(P_i\delta_i + \frac{1}{2}b_{ii}V_i^2 + Q_i\ln V_i).
\end{align}
Correspondingly,
\begin{align}
\label{L lineplusbus}
L = L_{\rm line} + L_{\rm bus},
\end{align}
where
\begin{align}
\label{laplacian line}
L_{\rm line}& =\frac{\partial^2 R_{\rm line}}{\partial z^2},\\
L_{\rm bus}& =\frac{\partial^2 R_{\rm bus}}{\partial z^2}.
\label{laplacian bus}
\end{align}
Note that $R_{\rm bus}$ contributes just in the diagonal terms of $L$
that are related with the algebraic variables $V_i$. Computing 
the differential of $L$ from (\ref{L lineplusbus}),
\begin{align}
\label{dL}
dL = dL_{\rm line} + dL_{\rm bus}.
\end{align}

In subsection \ref{Computing dLine} we compute $dL_{\rm line}$
in the new coordinates related to the transmission
lines and compute  $dL_{\rm bus}$ in the bus coordinates.

\subsection{Computing $dL_{\rm line}$ \label{Computing dLine}}

Similarly to (\ref{L and Ldash2}),
\begin{align}
L_{\rm line}  = H^TL'_{\rm line}H.
\end{align}
Then
\begin{align}
\nonumber
dL_{\rm line} & = d(H^TL'_{\rm line}H)\\
\label{dL line}
& =2H^TL'_{\rm line}dH+H^T(dL'_{\rm line})H. 
\end{align}

We first compute $L'_{\rm line}$.
$ R_{\rm line}$ can be nicely written  in the line coordinates as
\begin{align}
\label{Rdash lines coordenates}
 R'_{\rm line}&=-\sum_{k=1}^{\ell} b_{k} e^{\nu_k}\cos\theta_k.
 \end{align}
 Of course, since $ R_{\rm line}$ is a scalar, $ R'_{\rm line}=R_{\rm line}$. Then
\begin{align}
\label{laplacian dash}
L'_{\rm line}& =\frac{\partial^2 R'_{\rm line}}{\partial z'^2}
=
\left(
\begin{matrix}
\frac{\partial^2 R'_{\rm line}}{\partial \theta^2} & \frac{\partial^2 R'_{\rm line}}{\partial \theta \partial \nu}\\[5pt]
\frac{\partial^2 R'_{\rm line}}{\partial \nu \partial \theta} & \frac{\partial^2 R'_{\rm line}}{\partial \nu^2}
\end{matrix}
\right).
\end{align}

It is convenient to define for line $k$ the quantities 
\begin{align}
 \label{q}
 p_k&=\phantom{-}b_{k} e^{\nu_k}\sin\theta_k,\\
\label{p}
 q_k&=-b_{k} e^{\nu_k}\cos\theta_k.
\end{align}
$p_k$ is the real power flow on line $k$, and 
 $q_k$ is part of the expression for the reactive power flows on line $k$.  Then
\begin{align}
\label{L dash theta square}
\frac{\partial^2 R'_{\rm line}}{\partial \theta_k^2} & = b_{k} e^{\nu_k}\cos{\theta_k}=-q_{_k},\\
\label{L dash theta nu}
\frac{\partial^2 R'_{\rm line}}{\partial \theta_k\partial \nu_k} & = b_{k} e^{\nu_k}\sin{\theta_k}=p_{_k},\\
\label{L dash nu square}
\frac{\partial^2 R'_{\rm line}}{\partial \nu_k^2} & = -b_{k} e^{\nu_k}\cos{\theta_k}=q_{_k}.
\end{align}
And the block submatrices of $L'_{\rm line}$ in (\ref{laplacian dash}) become the diagonal matrices
\begin{align}
\label{L dash theta square2}
\frac{\partial^2 R'_{\rm line}}{\partial \theta^2} & =- {\rm diag} \{q_1,q_2,\ldots,q_{\ell}\},\\
\label{L dash theta nu2}
\frac{\partial^2 R'_{\rm line}}{\partial \theta\partial \nu} & = {\rm diag} \{p_1,p_2,\ldots,p_{\ell}\},\\
\label{L dash nu square2}
\frac{\partial^2 R'_{\rm line}}{\partial \nu^2} & =  {\rm diag} \{q_1,q_2,\ldots,q_{\ell}\}.
\end{align}
The corresponding block submatrices of $dL'_{\rm line}$ are then
\begin{align}
\label{dL dash theta square}
d\left(\frac{\partial^2 R'_{\rm line}}{\partial \theta^2} \right) & =- {\rm diag} \{dq_1,dq_2,\ldots,dq_{\ell}\},\\
\label{dL dash theta nu}
d\left(\frac{\partial^2 R'_{\rm line}}{\partial \theta\partial \nu}\right)  & = {\rm diag} \{dp_1,dp_2,\ldots,dp_{\ell}\},\\
\label{dL dash nu square}
d\left(\frac{\partial^2 R'_{\rm line}}{\partial \nu^2} \right) & =  {\rm diag} \{dq_1,dq_2,\ldots,dq_{\ell}\}.
\end{align}

Now we compute the matrix $dH$.
From (\ref{Hthetak}) and (\ref{Hnuk}) we have that
the entries different from zero of $dH$ are the entries related to (\ref{Hnuk}); i.e.,
\begin{align}
\label{dHnuk}
dH_{ki}={\displaystyle -\frac{dV_i}{V_i^2}}|A_{ik}|,\  i=m+1,\ldots,n, \ k=\ell+1,\ldots,2\ell.
\end{align}
Defining $V_i^{\ln}  =\ln{V_i}$, then
\begin{align}
\label{Vln}
dV_i^{\ln} & = d(\ln{V_i})=\frac{dV_i}{V_i}.
\end{align}
With (\ref{Vln}) in mind, (\ref{dHnuk}) becomes
\begin{align}
\label{dHnuk Vln}
dH_{ki}={\displaystyle -\frac{dV_i^{\ln}}{V_i}}|A_{ik}|,\  i=m+1,\ldots,n, \ k=\ell+1,\ldots,2\ell.
\end{align}

Lastly, we compute $dL_{\rm bus}$. Using (\ref{laplacian bus}), note that $L_{\rm bus}$ contributes only
in the diagonal terms of $L$ that are related with the algebraic 
variables $V_i$; i.e.,
\begin{align}
\label{Lbus}
(L_{\rm bus})_{ij}=\left\{
\begin{array}{ll}
{\displaystyle -b_{ii} + \frac{Q_i}{V^2_i}} & \mbox{if $i=j$ and $i=m+1,\ldots,n$}.\\[4pt]
 0 & \text{otherwise.}
\end{array} \right.
\end{align}
%
Using (\ref{Vln}),
\begin{align}
\label{dLbus}
(dL_{\rm bus})_{ij}=\left\{
\begin{array}{ll}
{\displaystyle -\frac{2Q_i}{V^2_i}dV_i^{\ln}} & \mbox{if $ i=j$ and  $i=m+1,\ldots,n$}.\\
 0 & \text{otherwise.}
\end{array} \right.
\end{align}

\subsection{Computing $x^TdLx$}

In this section we compute $x^TdLx$. From 
(\ref{dL}) we have
\begin{align}
\label{xTdlx} 
x^TdLx=x^T(dL_{\rm line})x + x^T(dL_{\rm bus})x.
\end{align}

First we calculate $x^T(dL_{\rm line})x$ using
(\ref{dL line}):
\begin{align}
\nonumber
x^T(dL_{\rm line})x & =x^T[H^T(dL'_{\rm line})H +2H^TL'_{\rm line}dH]x\\
\label{xTdLlinex0}
& = x'^{T}(dL'_{\rm line})x' +2x'^{T}L'_{\rm line}(dH)x.
\end{align}
Then
\begin{align}
\nonumber
x'^{T}&(dL'_{\rm line})x'  = (x'_{\theta},x'_{\nu})
\begin{pmatrix}
d\left(\frac{\partial^2 R_{\rm line}}{\partial \theta^2}\right) & d\left(\frac{\partial^2 R_{\rm line}}{\partial \theta \partial \nu}\right)\\[5pt]
d\left(\frac{\partial^2 R_{\rm line}}{\partial \theta \partial \nu}\right) & d\left(\frac{\partial^2 R_{\rm line}}{\partial \nu^2}\right)
\end{pmatrix}
\left(
\begin{matrix}
x'_{\theta} \\
x'_{\nu} 
\end{matrix}
\right),\\
\nonumber
&\mbox{and using (\ref{dL dash theta square}-\ref{dL dash nu square})},\\
\label{dL dash}
& =\sum_{k=1}^{\ell}[(x'_{\nu_k})^2 -(x'_{\theta_k})^2]dq_{_k}+
2\sum_{k=1}^{\ell}x'_{\theta_k}x'_{\nu_k}dp_{_k}.
\end{align}

Now we calculate $2x'^{T}L'_{\rm line}(dH)x$.

Using (\ref{dHnuk Vln}), the first $k$ entries of $(dH)x$ are zero, and the 
last $k$ entries of $(dH)x$  are 
\begin{align}
\label{dHx}
(dH)x=\begin{pmatrix}
0\\
c_{\nu}
\end{pmatrix},
\end{align}
where, writing $x_{_{V_{i}}}^{\ln}=\displaystyle\frac{x_{_{V_{i}}}}{V_i}$,
\begin{align}
\label{dHx2}
c_{\nu_k}
&=-\sum_{i=m+1}^{n}{\displaystyle x_{_{V_{i}}}^{\ln}dV_i^{\ln}}|A_{ik}|, \quad
k=\ell+1,\ldots,2\ell.
\end{align}
Now
\begin{align}
2x'^{T}(L'_{\rm line})(dH)x  &= 2 
(x'_{\theta},x'_{\nu})
\begin{pmatrix}
\frac{\partial^2 R_{\rm line}}{\partial \theta^2} &\frac{\partial^2 R_{\rm line}}{\partial \theta \partial \nu}\\[5pt]
\frac{\partial^2 R_{\rm line}}{\partial \theta \partial \nu}& \frac{\partial^2 R_{\rm line}}{\partial \nu^2}
\end{pmatrix}
\left(
\begin{matrix}
0 \\
c_{\nu}\end{matrix}
\right),\notag\\
\mbox{and using (\ref{L dash theta square}-\ref{L dash nu square})},\notag\\
\nonumber
&=2\sum_{k=1}^{\ell}\left(x'_{\theta_k}p_{_k}+ x'_{\nu_k}q_{_k}\right)c_{\nu_k}\\
\label{dHH}
&\hspace{-15mm}=-\!\sum_{i=m+1}^n\Bigg\{2\sum_{k=1}^{\ell}|A_{ik}|
(x'_{\theta_k}p_{_k} 
+ x'_{\nu_k}q_{_k})\,x_{_{V_{i}}}^{\ln}\Bigg\}dV_i^{\ln}. 
\end{align}
From (\ref{dL dash}) and (\ref{dHH}),
\begin{align}
\nonumber
 &x^T(dL_{\rm line})x  = 
 \sum_{k=1}^{\ell}[(x'_{\nu_k})^2 -(x'_{\theta_k})^2]dq_{_k}+
2\sum_{k=1}^{\ell}x'_{\theta_k}x'_{\nu_k}dp_{_k}\\
\label{xTdLlinex}
& \qquad -\sum_{i=m+1}^n\left\{2\sum_{k=1}^{\ell}|A_{ik}|
\left(x'_{\theta_k}p_{_k}+ x'_{\nu_k}q_{_k}\right)x_{_{V_{i}}}^{\ln}\right\}dV_i^{\ln}. 
\end{align}

Lastly, we compute $x^T(dL_{\rm bus})x$. From (\ref{dLbus}),
%
\begin{align}
\label{xTdLbusx} 
x^T(dL_{\rm bus})x= -\sum_{i=m+1}^n 2(x_{_{V_{i}}}^{\ln})^2Q_idV_i^{\ln}.
\end{align}

From (\ref{xTdLlinex}) and (\ref{xTdLbusx}),
\begin{align}
\nonumber
x^TdLx &= \sum_{k=1}^{\ell}[(x'_{\nu_k})^2 -(x'_{\theta_k})^2]dq_{_k}+
2\sum_{k=1}^{\ell}x'_{\theta_k}x'_{\nu_k}dp_{_k}\\
\nonumber
&\quad -\sum_{i=m+1}^n\left\{2\sum_{k=1}^{\ell}|A_{ik}|\left(x'_{\theta_k}p_{_k}+ x'_{\nu_k}q_{_k}\right)x_{_{V_{i}}}^{\ln}\right\}dV_i^{\ln}\\ 
\label{numeratordlambda}
& \quad -\sum_{i=m+1}^n2(x_{_{V_{i}}}^{\ln})^2Q_idV_i^{\ln}.
\end{align}

Expressing $dp_k$ and $dq_k$ in terms of $d\theta_k$, $dV_i^{\ln}$ and $dV_j^{\ln}$,
\begin{align}
\label{pthetaV}
dp_{_k} & = -q_{_k}d\theta_{k} + p_{_k}dV_i^{\ln} + p_{_k}dV_j^{\ln},\\
\label{qthetaV}
dq_{_k} & = p_{_k}d\theta_{k} + q_{_k}dV_i^{\ln} + q_{_k}dV_j^{\ln}.
\end{align}
Substituting (\ref{pthetaV}) and (\ref{qthetaV}) in (\ref{numeratordlambda}) and rearranging
terms,
\begin{align}
\nonumber
x^TdLx &= \sum_{k=1}^{\ell}[(x'_{\nu_k})^2 -(x'_{\theta_k})^2]
(p_{_k}d\theta_{k} + q_{_k}dV_i^{\ln} + q_{_k}dV_j^{\ln})\\
\nonumber
& \quad + 2\sum_{k=1}^{\ell}x'_{\theta_k}x'_{\nu_k}(-q_{_k}d\theta_{k} + p_{_k}dV_i^{\ln} + p_{_k}dV_j^{\ln})\\
\nonumber
&\quad -\sum_{i=m+1}^n\left\{2\sum_{k=1}^{\ell}|A_{ik}|\left(x'_{\theta_k}p_{_k}+ x'_{\nu_k}q_{_k}\right)x_{_{V_{i}}}^{\ln}\right\}dV_i^{\ln}\\ 
\nonumber
& \quad -\sum_{i=m+1}^n2(x_{_{V_{i}}}^{\ln})^2Q_idV_i^{\ln}\notag\\
\nonumber
&=\sum_{k=1}^{\ell}\left\{[(x'_{\nu_k})^2 -(x'_{\theta_k})^2]p_{_k} - 2x'_{\theta_k}x'_{\nu_k}q_{_k}\right\}d\theta_k\\
\nonumber
& \quad +\!\!\!\!\sum_{i=m+1}^n\sum_{k=1}^{\ell}|A_{ik}|\left\{x'_{\nu_k}[x'_{\nu_k}-2x_{_{V_{i}}}^{\ln}]-(x'_{\theta_k})^2\right\}q_{_k}dV_i^{\ln}\\
\nonumber
& \quad +\sum_{i=m+1}^n\sum_{k=1}^{\ell}|A_{ik}|\left\{x'_{\nu_k}-x_{_{V_{i}}}^{\ln}\right\}2x'_{\theta_k}p_{_k}dV_i^{\ln}\\
\label{numeratordlambdathetaV0}
&\quad -\sum_{i=m+1}^n2\left(x_{_{V_{i}}}^{\ln}\right)^2Q_idV_i^{\ln}.
\end{align}

Define
\begin{align}
\label{Cqk}
C_{q_{_{k}}} & = x'_{\nu_k}\left(x'_{\nu_k}-2x_{_{V_{i}}}^{\ln}\right) - (x'_{\theta_k})^2,\\
C_{p_{_{k}}} & = 2x'_{\theta_k}\left(x'_{\nu_k} - x_{_{V_{i}}}^{\ln}\right),\\
\label{Cpk}
C_{Q_i} & = -2\left(x_{_{V_{i}}}^{\ln}\right)^2.
\end{align}
Note that $C_{p_{_{k}}}\neq0$ only when the $k$th line is connecting two
load buses. Substituting (\ref{Cqk}-\ref{Cpk}) into
(\ref{numeratordlambdathetaV0}),
\begin{align}
\nonumber
&x^TdLx = \sum_{k=1}^{\ell}\left\{[(x'_{\nu_k})^2 -(x'_{\theta_k})^2]p_{_k} - 2x'_{\theta_k}x'_{\nu_k}q_{_k}\right\}d\theta_k\\
\label{numeratordlambdathetaV1}
& \quad + \sum_{i=m+1}^n\Bigg\{
\sum_{k=1}^{\ell}|A_{ik}|(C_{q_{_{k}}}q_{_k}+ C_{p_{_{k}}}p_{_k})+ C_{Q_i}Q_i\Bigg\}dV_i^{\ln}.
\end{align}
The expression (\ref{numeratordlambdathetaV1}) is the numerator of (\ref{dlambda}), so that
the final formula is 
\begin{align}
\nonumber
d\lambda &= -\frac{x^TdLx}{\alpha}  \notag\\&=-\frac{1}{\alpha}\left\{\sum_{k=1}^{\ell}
\left\{[(x'_{\nu_k})^2 -(x'_{\theta_k})^2]p_{_k} - 2x'_{\theta_k}x'_{\nu_k}q_{_k}\right\}d\theta_k \right.\notag\\
\label{dlambda partial theta partial V}
& \quad +\left. \sum_{i=m+1}^n\left[
\sum_{k=1}^{\ell}|A_{ik}|(C_{q_{_{k}}}q_{_k}+ C_{p_{_{k}}}p_{_k})+ C_{Q_i}Q_i \right]dV_{i}^{\ln} \right\},
\end{align}
where, repeating (\ref{alpha}) for convenience,
\begin{align}
\label{alphaagain}
\alpha=2\lambda x^TMx + x^TDx.
\end{align}
%
We make some general observations about formula (\ref{dlambda partial theta partial V}):
\begin{itemize}
\item Generator redispatch results in the changes $d\theta$ and $dV^{\rm ln}$ and affects only the numerator of (\ref{dlambda partial theta partial V}). The denominator $\alpha$ 
depends on the equivalent generator parameters, the eigenvalue $\lambda$, and the eigenvector $x$, and is the same 
for all generator redispatches. 
Thus, after accounting for the common effect of the denominator on all the generator 
redispatches,  we can determine the effective or ineffective generator dispatches by their varying effects on the numerator.
\item
The numerator depends on the redispatch via the changes in angles across the lines $d\theta$ and 
changes in load voltage magnitudes $dV^{\rm ln}$.
The coefficients of $d\theta$ depend on the mode eigenvector $x$ expressed in line coordinates and the 
real and reactive power line flows.
The coefficients of $dV^{\rm ln}$ additionally depend on the load voltage magnitude coordinates of the mode eigenvector $x$ 
and the reactive power load.
\item To find a good generator redispatch to improve the mode eigenvalue, we need to identify lines $k$  that have coefficients $d\theta_k$ of suitable magnitude and sign, and then find a redispatch that suitably changes  $\theta_k$ on those lines. 
If the redispatch also affects load voltage magnitudes, we need to also consider $dV^{\rm ln}$  and the coefficients of $dV^{\rm ln}$.
\item In (\ref{dlambda partial theta partial V}), the change in the mode eigenvalue is given as a complex number.
In practice, for maintaining oscillatory stability, we are most interested in the change in damping (the real part of $d\lambda$) or the change in the damping ratio.
\item The sensitivity of a mode depends linearly 
on the active and reactive power flow through every 
line of the network at the equilibrium. In
the case of load buses the constant reactive power
 demand also affects linearly the sensitivity
of the mode. 
\item The formula (\ref{dlambda partial theta partial V}) is independent of the scaling of the eigenvector $x$.
\end{itemize}
\section{Relating the redispatch to changes $d\theta$ and $dV_{i}^{\ln}$}
Formula (\ref{dlambda partial theta partial V}) relates the mode change $d\lambda$ to $d\theta$ and $dV_{i}^{\ln}$.
It remains to express $d\theta$ and $dV_{i}^{\ln}$ in terms of the redispatch $dP$ using the linearization 
of the load flow equations.

The linearization 
of the load flow equations is 
\begin{align}
\label{linearization load flow real}
 \sum_{j=1}^{2n-m} L_{ij}d z_j  & = dP_i, \quad i=1,2,\ldots,n,\\ 
\label{linearization load flow reactive}
\sum_{j=1}^{2n-m} L_{ij}d z_j  & = 0,\quad  i=m+1,\ldots,n.
\end{align}
where 
$d z=(d\delta,d V)^T$.
In matrix form,
\begin{align}
\label{linearization load flow matrix}
 L d z  & = \begin{pmatrix}
 dP\\0
 \end{pmatrix}
\end{align}
Then we can use the matrix pseudo-inverse (indicated by $\dagger$) to obtain
\begin{align}
\label{solve linearization load flow matrix}
\begin{pmatrix}
d\delta\\ dV
\end{pmatrix}=d z  & = L^{\dagger}
 \begin{pmatrix}
 dP\\0
 \end{pmatrix}
\end{align}
Then $d\theta$ and 
$dV^{\ln}$ are easily obtained:
\begin{align}
\label{final transform solve linearization load flow matrix}
d\theta_k&= \sum_{r=1}^n A_{rk}\, d\delta_r, \quad k=1,\ldots,\ell.\\
dV^{\ln}_i &=\frac{d V_i}{V_i},\quad i=m+1,\ldots,n.
\end{align}

\section{Computing $d\lambda$ for a 3-bus System}
\begin{figure}[h!]
\centerline{{\resizebox{8cm}{!}{%
\includegraphics{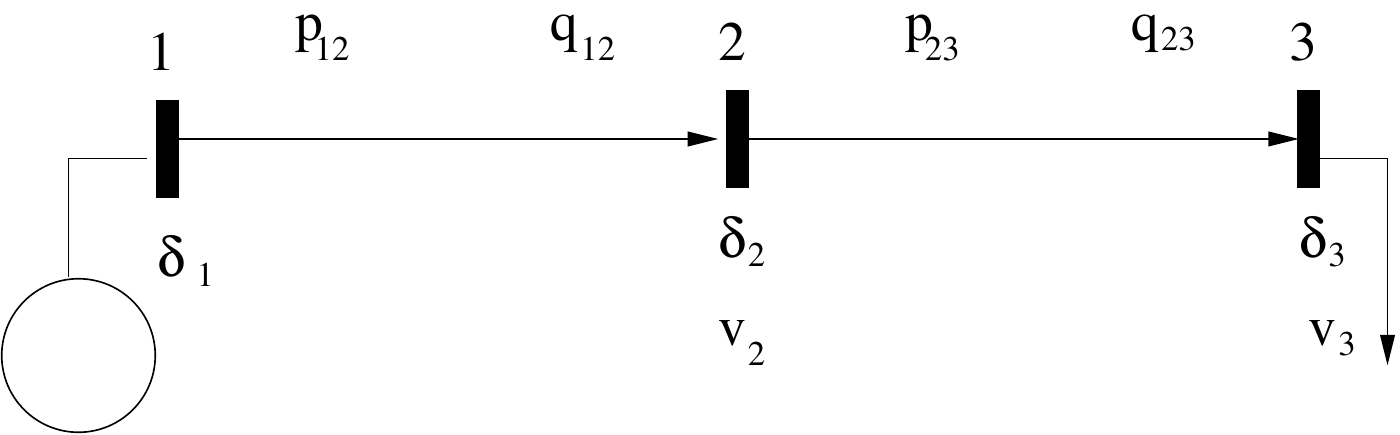}}}}
\caption{\small \label{3buses}3-bus system}
\end{figure}

In order to show an example of the derivation in a more explicit fashion,
we compute $d\lambda$ using
the new coordinates for the simple three bus system 
shown in Fig.~\ref{3buses}. Bus 1 is a generator bus,
bus 2 is a connecting point, and bus 3 is a load bus.

For this small system, $\delta=(\delta_1,\delta_2,\delta_3)^T$ and
$V=(V_2,V_3)^T$, so that  $z=(\delta_1,\delta_2,\delta_3,V_2,V_3)^T$.
The incidence matrix associated with the network is
\begin{align}
\label{A 3 buses}
A=\left(
\begin{array}{cc}
 1 & 0 \\
 -1 & 1 \\
 0 & -1 \\
\end{array}
\right),
\end{align}
so that the $z'=(\theta,\nu)$ coordinates are
\begin{align}
\theta=A^T\delta= \left(
\begin{array}{ccc}
 1 & -1 & 0 \\
 0 & 1 & -1 \\
\end{array}
\right)
 \left(
\begin{array}{c}
 \delta_1 \\
 \delta_2  \\
 \delta_3  \\
\end{array}
\right)
=
\left(
\begin{array}{c}
 \theta_1 \\
 \theta_2  \\
\end{array}
\right),
\end{align}
\begin{align}
\nu=|A^T| \left(
\begin{array}{c}
\ln{V_1} \\
\ln{V_2}  \\
\ln{V_3}  \\
\end{array}
\right)= \left(
\begin{array}{ccc}
 1 & 1 & 0 \\
 0 & 1 & 1 \\
\end{array}
\right)
 \left(
\begin{array}{c}
\ln{V_1} \\
\ln{V_2}  \\
\ln{V_3}  \\
\end{array}
\right)
=
\left(
\begin{array}{c}
\nu_1  \\
\nu_2  \\
\end{array}
\right).
\end{align}
Then the matrix transformation $H$ is
\begin{align}
\label{H 3 buses}
H=\frac{\partial h}{\partial z}
=\left(
\begin{array}{ccccc}
 1 & -1 & 0 & 0 & 0 \\
 0 & 1 & -1 & 0 & 0 \\
 0 & 0 & 0 & \frac{1}{V_2} & 0 \\
 0 & 0 & 0 & \frac{1}{V_2} & \frac{1}{V_3} \\
\end{array}
\right).
\end{align}
An eigenvector $\left(x_{\delta _1},x_{\delta _2},x_{\delta _3},x_{_{V_2}},x_{_{V_3}}\right)^T$
transforms as
\begin{align}
x'  =Hx&=\left(
\begin{array}{ccccc}
 1 & -1 & 0 & 0 & 0 \\
 0 & 1 & -1 & 0 & 0 \\
 0 & 0 & 0 & \frac{1}{V_2} & 0 \\
 0 & 0 & 0 & \frac{1}{V_2} & \frac{1}{V_3} \\
\end{array}
\right)
 \left(
\begin{array}{c}
x_{_{\delta _1}}\\
x_{_{\delta _2}}\\
x_{_{\delta _3}}\\
x_{_{_{V_2}}} \\
x_{_{_{V_3}}} \\
\end{array}
\right) \notag\\&=
\left(
\begin{array}{c}
 x_{_{\delta _1}}-x_{_{\delta _2}} \\
 x_{_{\delta _2}}-x_{_{\delta _3}} \\[3pt]
\displaystyle \frac{x_{_{_{V_2}}}}{V_2} \\
\displaystyle \frac{x_{_{V_2}}}{V_2}+\frac{x_{_{V_3}}}{V_3} \\
\end{array}
\right) =
\left(
\begin{array}{c}
x'_{\theta_1} \\
x'_{\theta_2} \\
x'_{\nu_1} \\
x'_{\nu_2} \\
\end{array}
\right).
\end{align}

The potential energy $R$ of the system is
\begin{align}
\label{R 3 buses0}
R & = R_{\rm line} + R_{\rm bus} =-\!\!\!\!\!\!\sum_{\substack{(i,j)=\\(1,2),(2,3)}}\!\!\!\!b_{ij}V_iV_j\cos(\delta_i - \delta_j)\notag\\
& \qquad - \sum_{i=1}^3 (P_i\delta_i + {\textstyle \frac{1}{2}}b_{ii}V_i^2 + Q_i\ln V_i). 
\end{align}
Expressing $R_{\rm line}$ in $z'$ coordinates we have
\begin{align}
\label{R 3 buses}
R & = R'_{\rm line} + R_{\rm bus} = -\sum_{k=1}^{2} b_{k} e^{\nu_k}\cos\theta_k\notag\\
& \qquad - \sum_{i=1}^3(P_i\delta_i + {\textstyle \frac{1}{2}}b_{ii}V_i^2 + Q_i\ln V_i). 
\end{align}

To compute the sensitivity of the mode, first we will get $x^TdLx$, the numerator
of (\ref{dlambda}), by using (\ref{xTdlx}), so
$x^TdLx=x^T(dL_{\rm line})x + x^T(dL_{\rm bus})x$. Working
with the first term, according to  (\ref{xTdLlinex0}) we have 
to compute $L'_{\rm line},~dL'_{\rm line}$ and $(dH)x$.

\begin{align}
&L'_{\rm line}=\notag\\
&\left(
\begin{array}{cccc}
b_{1} e^{\nu _1} \cos \theta _1 & 0 & b_{1} e^{\nu _1} \sin \theta _1 & 0 \\
 0 & b_{2} e^{\nu _2} \cos \theta _2& 0 & b_{2} e^{\nu _2} \sin \theta _2 \\
b_{1} e^{\nu _1} \sin \theta _1 & 0 & -b_{1} e^{\nu _1} \cos\theta _1 & 0 \\
 0 & b_{2} e^{\nu _2} \sin \theta _2 & 0 & -b_{2} e^{\nu _2} \cos \theta _2 \\
\end{array}
\right),\\
\nonumber
&\mbox{and, using (\ref{q}) and (\ref{p}),}\nonumber\\
\nonumber
& =
\left(
\begin{array}{cccc}
 -q_1 & 0 & p_1 & 0 \\
 0 & -q_2 & 0 & p_2 \\
p_1 & 0 & q_1 & 0 \\
 0 & p_2 & 0 & q_2 \\
\end{array}
\right).
\end{align}
Then
\begin{align}
 dL'_{\rm line}=
\left(
\begin{array}{cccc}
 -dq_1 & 0 & dp_1 & 0 \\
 0 & -dq_2 & 0 & dp_2 \\
dp_1 & 0 & dq_1 & 0 \\
 0 & dp_2 & 0 & dq_2 \\
\end{array}
\right)
\end{align}
and
\begin{align}
x'^{T}dL'_{\rm line}x' & =
\left(
\begin{array}{c}
x'_{\theta_1} \\
x'_{\theta_2} \\
x'_{\nu_1} \\
x'_{\nu_2} \\
\end{array}
\right)^T
\left(
\begin{array}{cccc}
 -dq_1 & 0 & dp_1 & 0 \\
 0 & -dq_2 & 0 & dp_2 \\
dp_1 & 0 & dq_1 & 0 \\
 0 & dp_2 & 0 & dq_2 \\
\end{array}
\right)
\left(
\begin{array}{c}
x'_{\theta_1} \\
x'_{\theta_2} \\
x'_{\nu_1} \\
x'_{\nu_2} \\
\end{array}
\right)\\
\label{x'TdLlinex' 3 buses}
& =
\sum_{k=1}^{2}[(x'_{\nu_k})^2 -(x'_{\theta_k})^2]dq_{_k}+
2\sum_{k=1}^{2}x'_{\theta_k}x'_{\nu_k}dp_{_k}.  
\end{align}

Now we  compute $2x'^{T}L(dH)x$
%
\begin{align}
dH= \left(
\begin{array}{ccccc}
 0 & 0 & 0 & 0 & 0 \\
 0 & 0 & 0 & 0 & 0 \\
 0 & 0 & 0 & -\frac{dV_2\vphantom{)}}{V_2^2\vphantom{)}} & \vphantom{\bigg)}0 \\
 0 & 0 & 0 & -\frac{dV_2}{V_2^2\vphantom{|}} & -\frac{dV_3}{V_3^2\vphantom{)}} \\
\end{array}
\right)
=
\left(
\begin{array}{ccccc}
 0 & 0 & 0 & 0 & 0 \\
 0 & 0 & 0 & 0 & 0 \\
 0 & 0 & 0 & -\frac{dV^{\ln}_2\vphantom{\big)}}{V_2\vphantom)} & 0 \\
 0 & 0 & 0 & -\frac{dV^{\ln}_2\vphantom{\big)}}{V_2} & -\frac{dV^{\ln}_3\vphantom{\big)}}{V_3} \\
\end{array}
\right),
\end{align}
Then $(dH)x$ is
\begin{align}
dHx & =
\left(
\begin{array}{ccccc}
 0 & 0 & 0 & 0 & 0 \\
 0 & 0 & 0 & 0 & 0 \\
 0 & 0 & 0 & -\frac{dV^{\ln}_2}{V_2} & 0 \\
 0 & 0 & 0 & -\frac{dV^{\ln}_2}{V_2} & -\frac{dV^{\ln}_3}{V_3} \\
\end{array}
\right)
 \left(
\begin{array}{c}
x_{_{\delta _1}}\\
x_{_{\delta _2}}\\
x_{_{\delta _3}}\\
x_{_{_{V_2}}} \\
x_{_{_{V_3}}} \\
\end{array}
\right)\notag\\& =
\left(
\begin{array}{c}
 0 \\
 0 \\
 -x^{\ln}_{_{V_2}}dV^{\ln}_2  \\
 - x^{\ln}_{_{V_2}}dV^{\ln}_2- x^{\ln}_{_{V_3}}dV^{\ln}_3 \\
\end{array}
\right).
\end{align}
And
\begin{align}
&2x'^{T}LdHx =  \notag\\
\notag
&2\left(
\begin{array}{c}
x'_{\theta_1} \\
x'_{\theta_2} \\
x'_{\nu_1} \\
x'_{\nu_2} \\
\end{array}
\right)^{\!\!\!T}
\left(
\begin{array}{cccc}
 -q_1 & 0 & p_1 & 0 \\
 0 & -q_2 & 0 & p_2 \\
p_1 & 0 & q_1 & 0 \\
 0 & p_2 & 0 & q_2 \\
\end{array}
\right)
\left(
\begin{array}{c}
 0 \\
 0 \\
 -x^{\ln}_{_{V_2}}dV^{\ln}_2  \\
 - x^{\ln}_{_{V_2}}dV^{\ln}_2- x^{\ln}_{_{V_3}}dV^{\ln}_3 \\
\end{array}
\right)
\\
& = 
-\!\sum_{i=2}^3\Bigg\{2\sum_{k=1}^{2}|A_{ik}|
(x'_{\theta_k}p_{_k} + x'_{\nu_k}q_{_k})(x_{_{V_{i}}}^{\ln})\Bigg\}dV_i^{\ln}. 
\label{2xdashLdHx}
\end{align}

From (\ref{x'TdLlinex' 3 buses}) and (\ref{2xdashLdHx}),
\begin{align}
\nonumber
x^{T}dL_{\rm line}x & =\sum_{k=1}^{2}[(x'_{\nu_k})^2 -(x'_{\theta_k})^2]dq_{_k}+
2\sum_{k=1}^{2}x'_{\theta_k}x'_{\nu_k}dp_{_k}\\
\label{xTdLinex 3 buses}
& \quad -\!\sum_{i=2}^3\Bigg\{2\sum_{k=1}^{2}|A_{ik}|
(x'_{\theta_k}p_{_k} + x'_{\nu_k}q_{_k})(x_{_{V_{i}}}^{\ln})\Bigg\}dV_i^{\ln}.
\end{align}

Now we compute $x^T(dL_{\rm bus})x$.

\begin{align}
L_{\rm bus}= \left(
\begin{array}{ccccc}
 0 & 0 & 0 & 0 & 0 \\
 0 & 0 & 0 & 0 & 0 \\
 0 & 0 & 0 & 0 & 0 \\
 0 & 0 & 0 & -b_{22} & 0 \\
 0 & 0 & 0 & 0 & -b_{33}+\frac{Q_3}{V_3^2} \\
\end{array}
\right),
\end{align}
and
\begin{align}
dL_{\rm bus}= \left(
\begin{array}{ccccc}
 0 & 0 & 0 & 0 & 0 \\
 0 & 0 & 0 & 0 & 0 \\
 0 & 0 & 0 & 0 & 0 \\
 0 & 0 & 0 & 0 & 0 \\
 0 & 0 & 0 & 0 & -2\frac{Q_3dV_3}{V_3^3} \\
\end{array}
\right)
=
\left(
\begin{array}{ccccc}
 0 & 0 & 0 & 0 & 0 \\
 0 & 0 & 0 & 0 & 0 \\
 0 & 0 & 0 & 0 & 0 \\
 0 & 0 & 0 & 0 & 0 \\
 0 & 0 & 0 & 0 & -2\frac{Q_3dV^{\ln}_3}{V_3^2} \\
\end{array}
\right).
\end{align}

Then
\begin{align}
x^TdL_{\rm bus}x & =  
\label{xTdLbusx3buses}
  -2\left(x_{_{V_{3}}}^{\ln}\right)^2Q_3\,dV_3^{\ln}.
\end{align}

From (\ref{xTdLinex 3 buses}) and (\ref{xTdLbusx3buses}),
\begin{align}
\nonumber
x^TdLx & = \sum_{k=1}^{2}[(x'_{\nu_k})^2 -(x'_{\theta_k})^2]dq_{_k}+
2\sum_{k=1}^{2}x'_{\theta_k}x'_{\nu_k}dp_{_k}\\
\nonumber
& \quad -\!\sum_{i=2}^3\Bigg\{2\sum_{k=1}^{2}|A_{ik}|
(x'_{\theta_k}p_{_k} + x'_{\nu_k}q_{_k})(x_{_{V_{i}}}^{\ln})\Bigg\}dV_i^{\ln}\\
\label{xTdLbusx3buses0}
& \quad 
-2\left(x_{_{V_{3}}}^{\ln}\right)^2Q_3dV_3^{\ln}.
\end{align}

From (\ref{pthetaV}-\ref{qthetaV}) and (\ref{Cqk}-\ref{Cpk}), (\ref{xTdLbusx3buses0}) 
becomes
\begin{align}
\nonumber
x^TdLx & = 
\sum_{k=1}^{2}\left\{[(x'_{\nu_k})^2 -(x'_{\theta_k})^2]p_{_k} - 2x'_{\theta_k}x'_{\nu_k}q_{_k}\right\}d\theta_k\\
& \quad + \sum_{i=2}^3\Bigg\{
\sum_{k=1}^{2}|A_{ik}|(C_{q_{_{k}}}q_{_k}+ C_{p_{_{k}}}p_{_k})+ C_{Q_i}Q_i\Bigg\}dV_i^{\ln}.
\end{align}
Then
\begin{align}
\nonumber
d\lambda & =-\frac{x^TdL x}{\alpha}\\
& \quad -\frac{1}{\alpha}\left\{
\sum_{k=1}^{2}\left\{[(x'_{\nu_k})^2 -(x'_{\theta_k})^2]p_{_k} - 2x'_{\theta_k}x'_{\nu_k}q_{_k}\right\}d\theta_k
\right.\\
& \quad + \sum_{i=2}^3\Bigg[
\left. \sum_{k=1}^{2}|A_{ik}|(C_{q_{_{k}}}q_{_k}+ C_{p_{_{k}}}p_{_k})+ C_{Q_i}Q_i\Bigg]dV_i^{\ln}\right\}.
\end{align}

\section{Special Case: Mode with Zero Damping}
\label{zerodamping}

To start to understand the general formula (\ref{dlambda partial theta partial V}), it is useful to consider special cases.
Given a pair $(\lambda,x)$ and following \cite{MalladaCDC11}, the quadratic equation
$\bar{x}^TQ(\lambda)x=0$ can be solved
to give 
\begin{align}
\label{quadratic2}
\lambda=\left\{
\begin{array}{ll}\displaystyle
-\frac{l(x)}{d(x)} & \text{if } m(x)=0, \\ \displaystyle
  \frac{-d(x)\pm\sqrt{d(x)^2-4m(x)l(x)}}{2m(x)} &  \text{otherwise, }
\end{array} \right.
\end{align}
\no where $m(x)=\bar{x}^TMx$, $d(x)=\bar{x}^TDx$ and
$l(x)=\bar{x}^TLx$. (\ref{quadratic2}) is the only calculation in the paper that makes use of 
$\bar{x}^TQ(\lambda)x=0$  instead of  $x^TQ(\lambda)x=0$.

Note that 
since $M \geq 0$, $D > 0$, and $L \geq 0$, we have 
$m(x) \geq 0$, $d(x) > 0$, and $l(x) \geq 0$ for all $ x$.

If the mode has zero damping so that  $\lambda$ is purely imaginary with  $\lambda= j\omega$, then from  (\ref{quadratic2}) 
we can see that $0=d(x)=\bar{x}^TDx=|\sqrt{D} x|^2$. 
Then $0=\sqrt{D} x=Dx$ and (\ref{Q}) becomes the real matrix equation 
$(-\omega^2 M + L)x =0$.
This implies that the  eigenvector $x$ (which is in general complex) can be taken to be real. Then the components 
of $x$ are  either exactly in phase or exactly $180\,^{\circ}$ out of phase 
with each other according to their sign.
Moreover, in this case, 
\begin{align}
{\alpha}&=j2\omega x^TMx,\\
\frac{-1}{\alpha}&=j\frac{1}{2\omega x^TMx},
\end{align}
and the formula for $d\lambda$ in (\ref{dlambda partial theta partial V}) becomes purely imaginary.
%
We conclude that in the case that
$\lambda$ is purely imaginary, changes in line angles and load bus voltage magnitudes  do 
not change the eigenvalue damping to first order; i.e, 
redispatch does neither stabilizes nor destabilizes  the operating
point. The only first order change possible in this case is a change in mode frequency.

This conclusion will remain approximately true if the mode 
damping is very small.
In the generic case of non-coincident eigenvalues,
since the eigenvector $x$ is a smooth function of parameters, 
it follows that a very lightly damped mode has an approximately real 
eigenvector $x$ and that the damping effect of redispatch is small.

\section{Special Case: Voltage Magnitudes Constant}

Another special case, for which the general formula (\ref{dlambda partial theta partial V}) simplifies dramatically,
is when the voltage magnitude is considered constant in all the buses.
The differential-algebraic
equations that describe the dynamics of the system are
(\ref{power balance}). Then (\ref{numeratordlambdathetaV0})
simplifies to %
\begin{align}
\label{active power case}
x^TdLx = -\sum_{k=1}^{\ell}(x'_{\theta_k})^2p_{_k}d\theta_k.
\end{align}
Substituting (\ref{active power case}) in (\ref{dlambda}), and letting $\lambda=\sigma \pm j\omega$, with $\omega$ positive,
\begin{align}
\label{sensitivity active power case}
d\lambda=d \sigma +jd \omega=
-\frac{x^TdLx}{\alpha} = \sum_{k=1}^{\ell}\frac{(x'_{\theta_k})^2p_{_k}}{\alpha}d\theta_k.
\end{align}
%
%
%
%
%

\subsection{Undamped mode case}

If the voltage magnitudes are assumed constant and $\lambda$ is a mode of the system with zero damping;
i.e., $\lambda=\pm j\omega$, then section~\ref{zerodamping} shows that 
the  eigenvector $x$ can be taken to be real and ${\alpha}=j2\omega x^TMx $.
Then (\ref{sensitivity active power case}) becomes
\begin{align}
\label{active power undamping0}
d \lambda=d \sigma +jd \omega
= -\sum_{k=1}^{\ell}j\frac{(x'_{\theta_k})^2p_{_k}}{2\omega m}d\theta_k
\end{align}
Since $(x'_{\theta_k})^2p_{_k}/(2\omega m)$ is a positive real number,
\begin{align}
\label{dsigma active power undamping}
d\sigma & = 0,\\
\label{domega active power undamping}
d\omega & =-\sum_{k=1}^{\ell}\left[\frac{(x'_{\theta_k})^2p_{_k}}{2\omega m}
\right]d\theta_k.
\end{align} 
In accordance with  section \ref{zerodamping}, 
 (\ref{dsigma active power undamping}) implies no change in $\sigma$ to first order. From (\ref{domega active power undamping}),
defining the positive number $(x'_{\theta_k})^2p_{_k}/2\omega m=a_k$ and 
substituting in (\ref{domega active power undamping}),
\begin{align}
\label{domega active power undamping1}
d\omega=-\sum_{k=1}^{\ell}a_k.d\theta_k=-a\cdotp d\theta.
\end{align}
Note that if $a$ and $d\theta$ are parallel, every entry of
the vector $a$ will contribute to $d\omega$. Which entries of 
the vector $a$ will contribute more? We answer this
question in subsections \ref{undamped mode 3 bus} and 
\ref{undamped mode n bus}.

\subsubsection{Undamped mode: 3-bus system \label{undamped mode 3 bus}}
\begin{figure}[h!]
\centerline{{\resizebox{8cm}{!}{%
\includegraphics{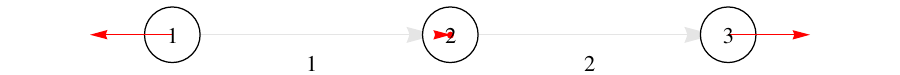}}}}
\caption{\small \label{3BusesSystem}The gray
lines  joining the buses show the magnitude of the  power flow  with the  grayscale 
and the direction of the  power flow 
with the arrows. Each line is numbered as shown. 
The red arrows at each bus show the oscillation mode shape associated with the critical eigenvalue
of the system; that is, the magnitude and direction of the 
real entries of the right eigenvector $x$ associated to
the critical eigenvalue $\lambda$. 1, 2 and 3 are generator buses.}
\end{figure}

In order to illustrate the use of formula (\ref{sensitivity active power case}),
we consider a simple 3-bus system with the power flow  and oscillating
mode pattern of its undamped critical mode shown in Fig.~\ref{3BusesSystem}.

The mode pattern shows that generator 1 is swinging against generator 3. Following the modal descriptions in 
\cite{FischerIREP}, 1 and 3 are antinodes of the system (locations with highest swing amplitude).
Generator 2 is not participating in the oscillation, so it is a node
(a location with zero swing amplitude). In more general power systems the nodes and antinodes may 
not be  located exactly at the buses.

 According to (\ref{domega active power undamping}), the 
sensitivity of the critical eigenvalue  is
\begin{align}
\label{active power undamping 3bus0}
d\omega= -\left[\frac{(x_{_{\delta_1}}-x_{_{\delta_2}})^2p_1}{2\omega m}\right]d\theta_1
-\left[\frac{(x_{_{\delta_2}}-x_{_{\delta_3}})^2p_2}{2\omega m}\right]d\theta_2,
\end{align}
where $p_1=b_{1}V_1V_2\sin{(\delta_1-\delta_2)}$ and 
$p_2=b_{2}V_2V_3\sin{(\delta_2-\delta_3)}$. As
the power flow goes from bus 1 to bus 2, $\delta_1>\delta_2$, then
$p_1>0$, and similarly $p_2>0$. As bus 2 is a node, $x_{_{\delta_2}}=0$, and
\begin{align}
\label{active power undamping 3bus}
d\omega= -\left[\frac{x_{_{\delta_1}}^2p_1}{2\omega m}\right]d\theta_1
-\left[\frac{x_{_{\delta_3}}^2p_2}{2\omega m}\right]d\theta_2.
\end{align}
Defining the positive real numbers $a_1=x_{_{\delta_1}}^2p_1/(2\omega m)$
and $a_2=x_{_{\delta_3}}^2p_2/(2\omega m)$ and substituting in
(\ref{active power undamping 3bus}),
\begin{align}
\label{domega 3bus}
 d\omega=-a_1d\theta_1-a_2d\theta_2=-a\cdotp d\theta,
\end{align}
where $a=(a_1,a_2)^T$, $d\theta=(d\theta_1,d\theta_2)^T$.
Define $\omega_{i}$ as the natural frequency of the system
in the base case; i.e., in the case of zero redispatch. 
Define $\omega_{f}$ as the natural frequency of the system
after redispatch, so that $d\omega=\omega_f-\omega_i$ Then 
$\omega_f=\omega_i+d\omega$. There are several cases:
\begin{enumerate}
\item Transfer between an antinode and an  antinode. There are
two subcases:
\begin{enumerate}
\item The transfer is made in the direction of the
power flow in the base case; i.e., from bus 1 to 
bus 3. Then the vectors $a$ and $d\theta$ are parallel. And
from (\ref{domega 3bus}), $d\omega<0$ and $\omega_f<\omega_i$, 
so the frequency of the mode decreases with the redispatch.
\item The transfer is made in the opposite direction
of the power flow in the base case; i.e., from bus
3 to bus 1. Then $a$ and $d\theta$ are antiparallel.
From (\ref{domega 3bus}), $d\omega>0$ and $\omega_f>\omega_i$, 
so the frequency of the mode increases with the redispatch.
\end{enumerate}
\item Transfer between a  node and an antinode;
for example, between bus 1 and 2. From (\ref{domega 3bus}), 
if the transfer is made in the direction of the base case power flow, then  $d\theta$ is positive and $\omega_f<\omega_i$.
If the transfer is made in the opposite direction to the power flow
in the base case, then $d\theta<0$ and $\omega_f>\omega_i$; i.e, 
the frequency increases with the redispatch.
\end{enumerate}
From  cases 1 and 2 we can conclude that the frequency of the 
mode decreases when the vectors $a$ and $d\theta$ are 
parallel. As $a$ is a vector with positive real entries,
to decrease the frequency the redispatch has to be done in the 
same direction as the power flow in the base case.

%

\subsubsection{Undamped mode: n-bus system  \label{undamped mode n bus} }

We consider an n-bus system that has an interarea
mode with zero damping; i.e, $\lambda=\pm j\omega$.
Then
\begin{align}
\label{puraly imaginary}
d\omega = -\sum_{k=1}^{\ell}\left[\frac{(x'_{\theta_k})^2p_{_k}}{2\omega m(x)}
\right]d\theta_k=-\sum_{k=1}^{\ell}a_k d\theta_k.
\end{align}
We note that $a_k \geq 0~\mbox{with}~k=1,\ldots \ell$.
If vectors $a$ and $d\theta$, are parallel (i.e., every
$d\theta_k >0$, or, in other words,  the redispatch causes power in every line 
to increase
in the direction of the power flow in the base case), then
every entry of the summation in (\ref{puraly imaginary}) will contribute to 
the decrease of the frequency of the mode. 
Any lines for which the redispatch causes the power 
to decrease
in the direction of the power flow in the base case will tend to 
increase the frequency of the mode.

The terms of
the summation (\ref{puraly imaginary}) that contribute more
correspond to those lines in which the product
$(x'_{\theta_k})^2p_{_k}$ is large.
These lines have large power flows and a large change in the eigenvector angle across the line.

One case of interest is when there is a power system area that includes an antinode $A_1$
transferring power to another power system area that includes an antinode $A_2$, but $A_2$ is swinging in 
the opposite direction to $A_1$.
Consider a path of lines joining $A_1$ to $A_2$ in which the power flow in each line is in the direction 
from $A_1$ to $A_2$.  Also assume that the amplitude of the oscillation behaves sinusoidally in space so that it decreases as one moves
on the path away from antinode $A_1$ until a node $N$ is encountered, and then the amplitude increases,
but with opposite phase as one passes from the node $N$ to antinode $A_2$.
Since antinodes are maxima of oscillation amplitude, near the antinode, changes in the eigenvector components are 
small and $(x'_{\theta_k})^2$ is small.
At the node the amplitude of the oscillation is zero but the gradient of the change in amplitude is large, and $(x'_{\theta_k})^2$ is large.
Thus if there is redispatch from $A_1$ to $A_2$ that increases the power flow in all the lines in the path, 
then the lines in the path near node $N$ contribute the most to decreasing the frequency of the mode.
A redispatch from $A_1$ to $N$, or a redispatch from $N$ to $A_2$ will also decrease the frequency of the mode.

%
%
%

\subsection{Damped mode case \label{damped mode case}}

Interarea modes are lightly damped electromechanical modes
of oscillation. In this section the sensitivity of a lightly
damped mode will be treated.
The sensitivity of a mode is given by (\ref{sensitivity active power case}).
We write $\alpha$ as
\begin{align}
\label{alpharandI}
\alpha = \alpha_r + j \alpha_{_I}.
\end{align}
Substituting (\ref{alpharandI}) in (\ref{sensitivity active power case}),
\begin{align}
d\lambda&=d\sigma +jd\omega=
\left[\frac{\alpha_r-j\alpha_{_I}}{\alpha_r^2+\alpha_{_I}^2}\right]\sum_{k=1}^{\ell}(x'_{\theta_k})^2p_{_k}d\theta_k \\
& =\left[\frac{\alpha_r-j\alpha_{_I}}{\alpha_r^2+\alpha_{_I}^2}\right]\sum_{k=1}^{\ell}(Re[(x'_{\theta_k})^2]+jIm[(x'_{\theta_k})^2])p_{_k}\\
& = \sum_{k=1}^{\ell}\frac{\alpha_rRe[(x'_{\theta_k})^2]+\alpha_{_I}Im[(x'_{\theta_k})^2]}{\alpha_r^2+\alpha_{_I}^2}p_{_k}\notag\\
& \qquad +j\sum_{k=1}^{\ell} \frac{\alpha_rIm[(x'_{\theta_k})^2]-\alpha_{_I}Re[(x'_{\theta_k})^2]}{\alpha_r^2+\alpha_{_I}^2}p_{_k}.
\end{align}
Then
\begin{align}
\label{dsigma active power damping}
d\sigma & = \sum_{k=1}^{\ell} 
\frac{\alpha_{_I}Im[(x'_{\theta_k})^2]+\alpha_rRe[(x'_{\theta_k})^2]}{\alpha_r^2+\alpha_{_I}^2}p_{_k}d\theta_k\notag\\&
= \sum_{k=1}^{\ell} a_{rk}d\theta_k=a_{_r}\cdotp d\theta,\\
\label{domega active power damping}
d\omega & = \sum_{k=1}^{\ell}
\frac{\alpha_rIm[(x'_{\theta_k})^2]-\alpha_{_I}Re[(x'_{\theta_k})^2]}{\alpha_r^2+\alpha_{_I}^2}p_{_k}d\theta_k\notag\\&
= \sum_{k=1}^{\ell} a_{_{Ik}}d\theta_k=a_{_I}\cdotp d\theta.
\end{align}
The ideal case to increase the magnitude of $\sigma$ and decrease $\omega$
(and with this increase the damping ratio) is when $a_r$ and $a_{_I}$ 
are parallel vectors and antiparallel with the vector $d\theta$. 
If $d\theta$ is antiparallel just with $a_r$, $\sigma$ will 
increase, but also $\omega$ will increase which is not good.
If $d\theta$ is antiparallel just with $a_{_I}$, $\omega$ will
decrease, but also $\sigma$ will decrease which is also not good.
Which entries of the vectors $a_r$ 
and $a_{_I}$ will contribute more?. We answer this
question in subsections \ref{damped mode 3 bus} and 
\ref{damped mode n bus}.

\subsubsection{Damped mode: 3-bus system \label{damped mode 3 bus}}
\begin{figure}[h!]
\centerline{{\resizebox{8cm}{!}{%
\includegraphics{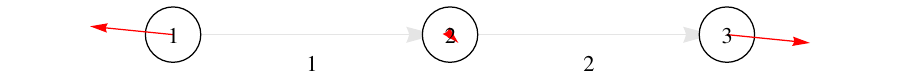}}}}
\caption{\small \label{3BusDamped}The gray
lines  joining the buses show the magnitude of the  power flow  with the  grayscale 
and the direction of the  power flow 
with the arrows. Each line is numbered as shown. 
The red arrows at each bus show the oscillation mode shape; 
that is, the magnitude and direction of the 
complex entries of the right eigenvector $x$ associated to
the critical complex eigenvalue $\lambda$. Buses 1, 2 and 3 are generator buses.}
\end{figure}

In this section the sensitivity of the lightly damped
electromechanical mode of oscillation of a 3-bus system 
is treated. The power flow and oscillating mode
pattern of its critical mode is shown in 
Fig.~\ref{3BusDamped}.

The mode pattern shows that generator 1 is swinging against generator 3 and that
bus 2 is not participating in the oscillation.
According to (\ref{dsigma active power damping})
and (\ref{domega active power damping}), the sensitivity of the nonzero eigenvalue 
of the system is given by
\begin{align}
\label{dsigma active power damping 3 buses0}
d\sigma & =
\frac{\alpha_{_I}Im[(x'_{\theta_1})^2]+\alpha_rRe[(x'_{\theta_1})^2]}{\alpha_r^2+\alpha_{_I}^2}p_{_1}d\theta_1\notag\\
&\qquad + \frac{\alpha_{_I}Im[(x'_{\theta_2})^2]+\alpha_rRe[(x'_{\theta_2})^2]}{\alpha_r^2+\alpha_{_I}^2}p_{_2}d\theta_2\\&
=a_rd\theta,\\
\label{domega active power damping 3 buses0}
d\omega & = \frac{\alpha_rIm[(x'_{\theta_1})^2]-\alpha_{_I}Re[(x'_{\theta_1})^2]}{\alpha_r^2+\alpha_{_I}^2}p_{_1}d\theta_1\notag\\
& \qquad +\frac{\alpha_rIm[(x'_{\theta_2})^2]-\alpha_{_I}Re[(x'_{\theta_2})^2]}{\alpha_r^2+\alpha_{_I}^2}p_{_2}d\theta_1\\&
=a_{_I}d\theta,
\end{align}
where $p_1=b_{1}V_1V_2\sin{(\delta_1-\delta_2)}$ and 
$p_2=b_{2}V_2V_3\sin{(\delta_2-\delta_3)}$. As
the power flow goes from bus 1 to bus 2, $\delta_1>\delta_2$ so that
$p_1>0$. Similarly, $p_2>0$.

From Fig.~\ref{3BusDamped} we can see that $x_{\delta_1}$ is
in the second quadrant of the complex plane and that
$x_{\delta_3}$ is in the fourth quadrant of the complex plane. Then

\begin{enumerate}
\item The complex numbers $x^TMx,~x^TDx$ are in the 
fourth quadrant of the complex plane. Then
\begin{align}
\label{alpha r and I}
\alpha &=2\lambda x^TMx + x^TDx\notag\\&=2(-\sigma+j\omega)x^TMx + x^TDx\notag\\&
=\alpha_r + j \alpha_{_I},
\end{align}
with $\alpha_r,\alpha_{_{I}}$ positive real numbers 
and $\alpha_r \ll \alpha_{_I}$.
\item $a_{_{I1}}<0$, $a_{_{I2}}<0$. So
from (\ref{domega active power damping 3 buses0}) to decrease $\omega$,
the redispatch has to be done in the direction of the power flow
in the base case. This result coincides with the conclusions for the 
undamped mode case.
\item $Re[(x'_{\theta_1})^2]>0$, $Re[(x'_{\theta_2})^2]>0$, 
$Im[(x'_{\theta_1})^2]<0$, $Im[(x'_{\theta_2})^2]<0$, so to
increase $|\sigma|$ we have to make the redispatch through the line 
in which the entry of $a_r$ is negative.
\end{enumerate}
Note that $|d\sigma|<|d\omega|$.


\subsubsection{Damped mode: $n$-bus system  \label{damped mode n bus}}

The sensitivity of an electromechanical mode
of oscillation of a network of $n$ buses is
given by equations (\ref{dsigma active power damping}) and
(\ref{domega active power damping}). The ideal case to increase 
the magnitude of $\sigma$ and decrease $\omega$ 
(and with this increase the damping ratio) is when $a_r$, $a_{_I}$ 
are parallel vectors and antiparallel with the vector $d\theta$.
If $d\theta$ is antiparallel just with $a_r$, $\sigma$ will 
increase, but also $\omega$ will increase which is not good.
If $d\theta$ is antiparallel just with $a_{_I}$, $\omega$ will
decrease, but also $\sigma$ will decrease which is also not good.
The terms of the summations (\ref{dsigma active power damping}) and
(\ref{domega active power damping}) that contribute more are those in which the product
$(x'_{\theta_k})^2p_{_k}$ is large. 
We would expect, as discussed in section \ref{undamped mode n bus}, that $(x'_{\theta_k})^2p_{_k}$ would be large 
in lines with substantial power flows that are near nodes at which the oscillation phase changes by approximately 180 degrees.
The redispatch should be chosen to exploit these lines, but we need to learn more about the general spatial 
structure of the modes to be able to better describe  this with confidence  and in detail.
%
%

\section{Verifying the new formula: AC power flow, 10-bus system}
\label{verify}

\begin{figure}[!h]
\includegraphics[width=\columnwidth]{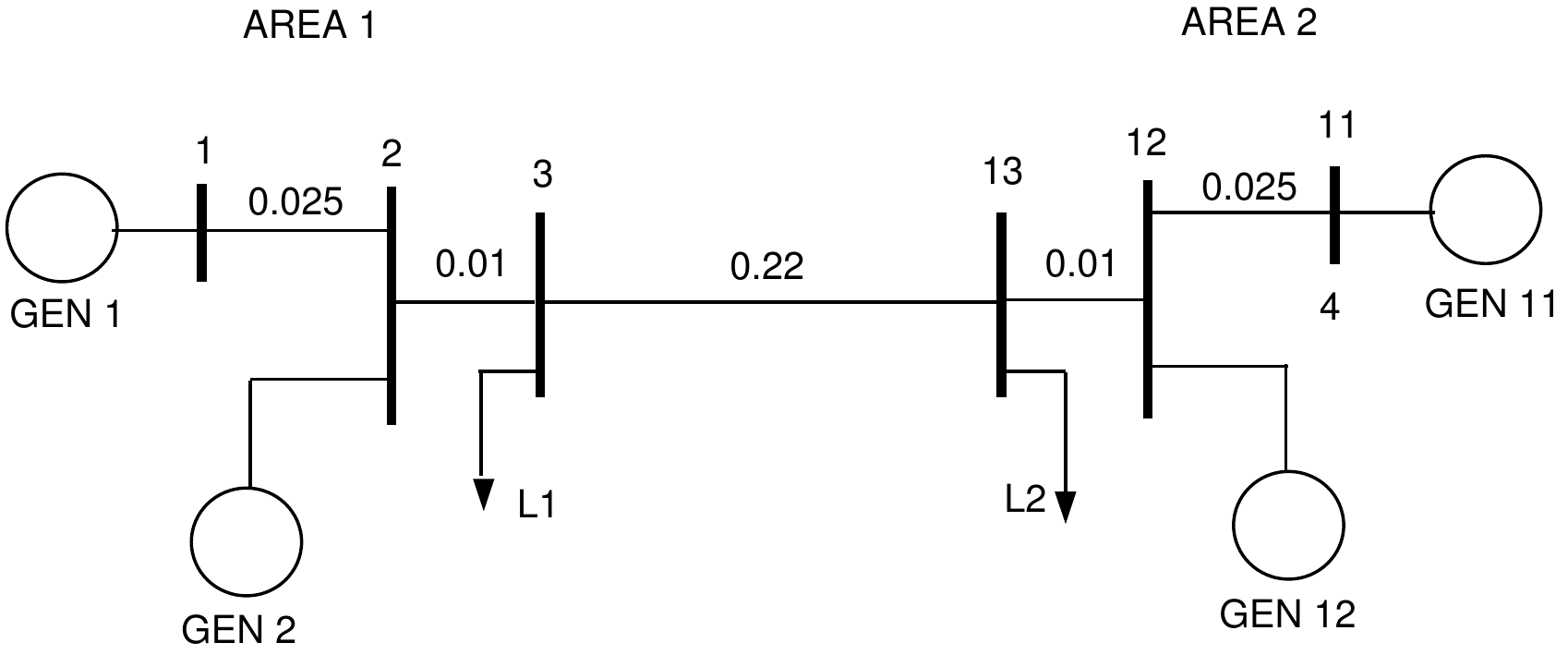}
\caption{\small \label{10BusFigure}10-bus system}
\end{figure}

In this section,  formula (\ref{dlambda partial theta partial V}) 
is verified in the 10-bus system shown in  Fig.~\ref{10BusFigure}. 
The system is based on the system in \cite{Klein91}, and consists of two similar areas connected by a weak tie line.
Each generator is represented by the same 
classical model with $H=6$.5 s, $D=1$.0 s, and transient reactance
$x'=0$.3. The internal constant voltage magnitudes of the 
generators are $V_1=0$.998337, $V_2=1$.26781, $V_3=1$.0782
and $V_4=1$.1449.
In the base case, $p_{_7}=3$.8897 is flowing 
through the tie line from area 1 to area 2. Table 
\ref{Generation Table} shows the generation and the power demanded by 
the constant loads in the base case.

\begin{table}[h]
\caption{\label{Generation Table}{Generator and load bus data of 
10-bus system
} }
\centering \begin{tabular}{cccccc}
bus & type & $P_g$ &  $P_L$ & $Q_L$ \\
\hline
 1 & G & 7.0  & 0.0 & 0.0\\
 2 & G & 7.0  & 0.0 & 0.0\\
 3 & G & 7.22049 & 0.0 & 0.0 \\
 4 & G & 7.0 & 0.0 & 0.0\\
 5 & L & 0.0 & 10.110245 & 1.0\\
 6 & L & 0.0 & 18.110245 & 1.0\\
 \hline
\end{tabular}
\end{table}

All the numerical computation is done with  the 
software Mathematica. First the power flow equations
are solved, and then the base case eigenvalues are computed. The system has three electromechanical
modes. Table \ref{Eigenvalues 10 Buses Base Case} shows the
electromechanical eigenvalues of the system for the base case.

\begin{table}[!h]
\caption{\label{Eigenvalues 10 Buses Base Case}{Eigenvalues of 10-bus system
in the base case
} }
\centering \begin{tabular}{ccl}\\
mode base case&
eigenvalue (rad/s)& Swing profile
\\
\hline
$\lambda_{1i}$ & -0.038462 + 8.8206i & 1,4 $\leftrightarrow$ 2,3 \\
$\lambda_{2i}$ & -0.038462 + 8.6023i & 1,4 $\leftrightarrow$ 2,3 \\
$\lambda_{3i}$ & -0.038462 + 2.3832i & 1,2 $\leftrightarrow$ 3,4 \\
\hline
\end{tabular}
\end{table}

\begin{figure}[!h]
\includegraphics[width=\columnwidth]{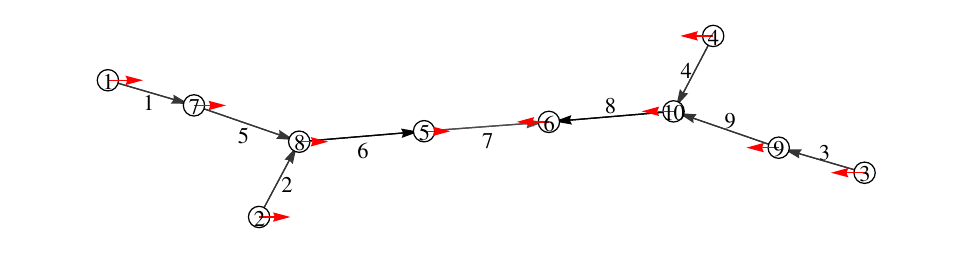}
\caption{\small \label{10BusesActivePowerThroughLines}The gray
lines  joining the buses show the magnitude of the  power flow  with the  grayscale 
and the direction of the  power flow 
with the arrows. Each line is numbered as shown. 
The red arrows at each bus show the oscillation mode shape; that is, the magnitude and direction of the 
 entries of the right eigenvector $x_{_{\delta}}$ associated with
the complex eigenvalue $\lambda_{3i}$. 1, 2, 3 and 4 are generator buses and
5 and 6 are load buses.}
\end{figure}

The power flow and oscillation for the base case
is shown in Fig.~\ref{10BusesActivePowerThroughLines}
as well as  the mode pattern of $\lambda_{3i}$. The mode pattern
shows that area 1 is swinging against area 2.

\begin{table}[h!]
\caption{\label{Eigenvalue lambda3 10 buses}{$\lambda_{3f}$ for redispatch
from G1 to G3 in 10-bus system
} }
\centering\begin{tabular}{ccc}\\
\!\!\!\!Redispatch\!\!\!\!\!\!\!\! & Exact mode & Approximate mode
\\
\hline
 0.000 & -0.038462 + 2.3832j & -0.038462 + 2.3832j \\
 0.003 &-0.038462 + 2.3785j & -0.038462 + 2.3786j\\
 0.006 & -0.038462 + 2.3738j & -0.038462 + 2.3739j\\
 0.009 & -0.038462 + 2.3691j & -0.038462 + 2.3692j\\
 0.010 & -0.038462 + 2.3675j & -0.038462 + 2.3676j\\
 0.03 & -0.038462 + 2.3350j & -0.038462 + 2.3357j\\
 0.06 & -0.038462 + 2.2829j &-0.038462 + 2.2858j\\
 0.09 & -0.038462 + 2.2262j & -0.038462 + 2.2331j \\
 0.10 & -0.038462 + 2.2061j & -0.038462 + 2.2149j\\
 0.15 & -0.038462 + 2.0947j &-0.038462 + 2.1173j \\
 0.20 & -0.038462 + 1.9586j & -0.038462 + 2.0060j\\
 0.25 & -0.038462 + 1.7810j &-0.038462 + 1.8735j \\
 0.30 & -0.038462 + 1.5152j & -0.038462 + 1.7005j\\
 \hline
\end{tabular}
\end{table}

\begin{figure}[!h]
\includegraphics[width=\columnwidth]{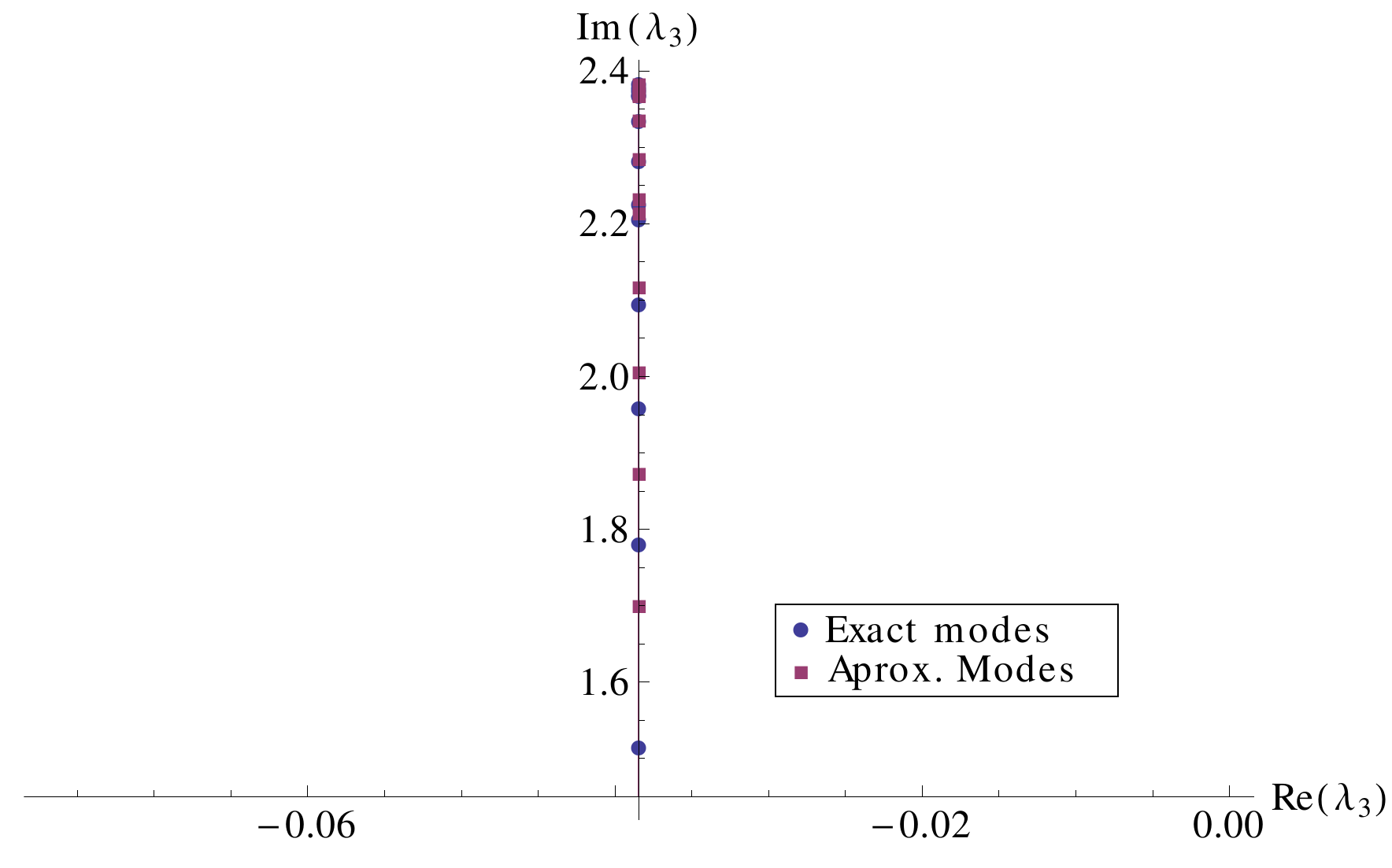}
\caption{\small \label{10BusesEigenvalues} Comparing the
exact and approximate modes in  the 10-bus system.}
\end{figure}

\begin{figure}[!h]
\includegraphics[width=\columnwidth]{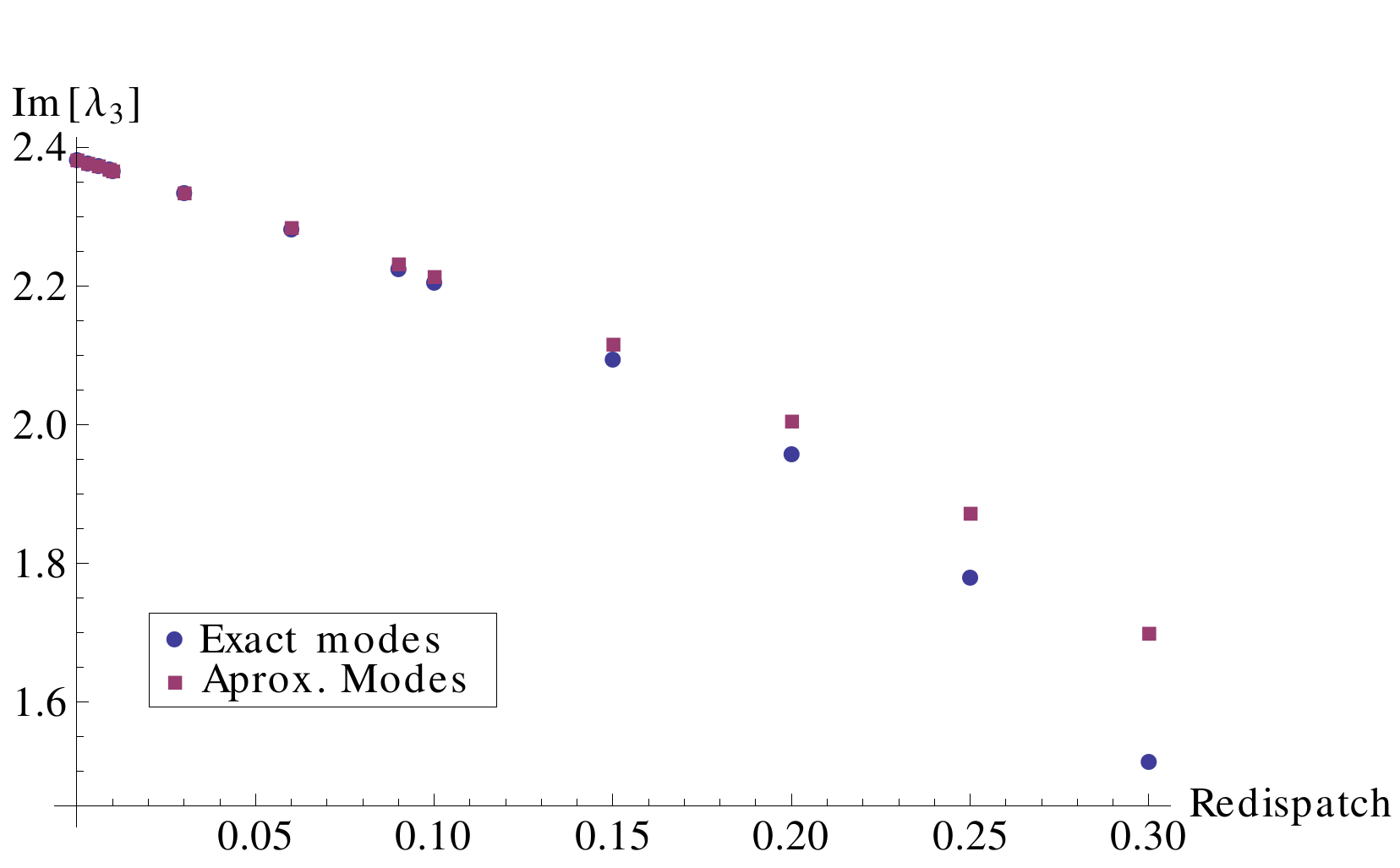}
\caption{\small \label{10BusesEigenvaluesRedispatch} Exact and approximate
mode frequencies versus amount of redispatch in the 10-bus system.}
\end{figure}

We examine changes in $\lambda_{3i}$  to test formula (\ref{dlambda partial theta partial V}). 
Redispatch is made between generator 1 of area 1 and generator 3 of
area 2. The generation of G1 is increased by an amount $r$ and the generation of 
G3 is decreased by  $r$. Using formula (\ref{dlambda partial theta partial V}),
$d\lambda_{3}$ is computed for several values of $r$, then the approximate eigenvalue
$\lambda_{3f}=d\lambda_{3} + \lambda_{3i}$ is calculated for every $r$. 
Table \ref{Eigenvalue lambda3 10 buses} shows
$\lambda_{3f}$ for different steps of redispatch between G1 and G3 and compares the 
exact and approximate eigenvalues.
Fig.~\ref{10BusesEigenvalues} compares the exact and approximate eigenvalues
of table \ref{Eigenvalue lambda3 10 buses} in the complex plane
and Fig.~\ref{10BusesEigenvaluesRedispatch} compares
the exact and approximate imaginary part of the eigenvalues versus 
the redispatch. From table \ref{Eigenvalue lambda3 10 buses} we
can confirm that formula (\ref{dlambda partial theta partial V}) reproduces the first order variation of 
the eigenvalues with respect to the redispatch. 

\section{6-bus system \label{6 bus system}}
\looseness=-1

In this section, we illustrate the use
of formulas (\ref{dsigma active power damping}) and
(\ref{domega active power damping}) in a simple 
6-bus system. These formulas compute the sensitivity
in the special case in which the voltage
is considered constant at every bus. The loads 
are modeled with frequency dependence of real 
power. The bus data of the  system is given 
in the table \ref{Generation Table 6 Bus} and the 
data of the transmission lines is given in table 
\ref{Transmission Table 6 Bus}.

\begin{table}[h!]
\caption{\label{Generation Table 6 Bus}{Bus data
of the 6-bus system
} }
\centering \begin{tabular}{cccccc}\\
bus & type & H (s) & D (s) &$P_g$ & $P_L$  \\
\hline
 1 & G & 3.0 & 2.0 & 0.8  & 0.0 \\
 2 & G & 3.0 & 2.0 & 0.8  & 0.0 \\
 3 & G & 24.0 & 16.0 & 6.4  & 0.0 \\
 4 & L & 0.0 & 2.0 & 0.0 & 1.0 \\
 5 & L & 0.0 & 2.0 & 0.0 & 1.0 \\
 6 & L & 0.0 & 16.0 & 0.0 & 6.0 \\
 \hline
\end{tabular}
\end{table}

\begin{table}[h!]
\caption{\label{Transmission Table 6 Bus}{Transmission line data
of the 6-bus system
} }
\centering \begin{tabular}{cl}
Line & ~~$x$ \\
\hline
 1 & 0.45  \\
 2 & 0.45 \\
 3 & 0.0563 \\
 4 & 0.02 \\
 5 & 0.075 \\
 \hline
\end{tabular}
\end{table}

The system has two electromechanical modes. 
Table \ref{Eigenvalues 6 Buses Base Case} shows the
electromechanical eigenvalues of the system for the base case.

\begin{table}[!h]
\caption{\label{Eigenvalues 6 Buses Base Case}{Eigenvalues of the 6-bus system
in the base case
} }
\centering \begin{tabular}{ccccc}
&&&&swing\\[-2pt]
 & f (Hz) & $\zeta (\%)$ & eigenvalue (rad/s) & profile
\\
\hline
$\lambda_{1i}$ & 1.53802 & 1.81694 & -0.175611 + 9.66364j & 1,2$\leftrightarrow$3\\
$\lambda_{2i}$ & 1.72281 & 1.54097 & -0.166826 + 10.8247j & 1 $\leftrightarrow$ 2 
 \\
\hline
\end{tabular}
\end{table}

\begin{figure}[h!]
\centerline{{\resizebox{8cm}{!}{%
\includegraphics{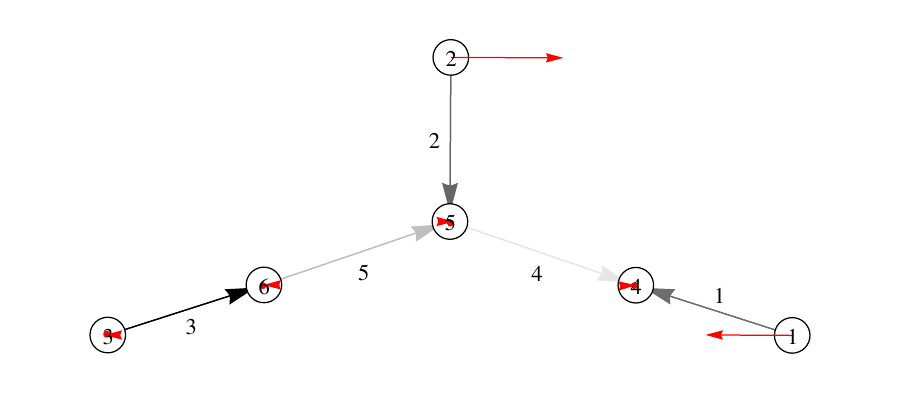}}}}
\caption{\small \label{ExpositoPowerThroughLines}Six-bus system: 
The gray
lines  joining the buses show the magnitude of the  power flow  with the  grayscale and the direction of the  power flow 
with the arrows. Each line is numbered as shown. 
The red arrows at each bus show the oscillation mode shape; that is,  the magnitude and direction of the 
complex entries of the right eigenvector $x_{_\delta}$ associated to
the critical complex eigenvalue $\lambda_{2i}$. Buses 1 and 2 are antinodes and
buses 3,4,5,6 are nodes. }
\end{figure}

The power flow oscillation in the base case and the mode pattern 
of $\lambda_{2i}$ are shown in Fig. \ref{ExpositoPowerThroughLines}. 
The mode pattern shows that G1 is swinging against G2.
The real coefficients $a_{_{rk}}$ and $a_{_{Ik}}$
in the equations (\ref{dsigma active power damping}) and
(\ref{domega active power damping}) for the 6-bus system are shown in table 
\ref{ar and aI 6 Buses Base Case}.

\begin{figure}[h!]
\small
\includegraphics[width=\columnwidth]{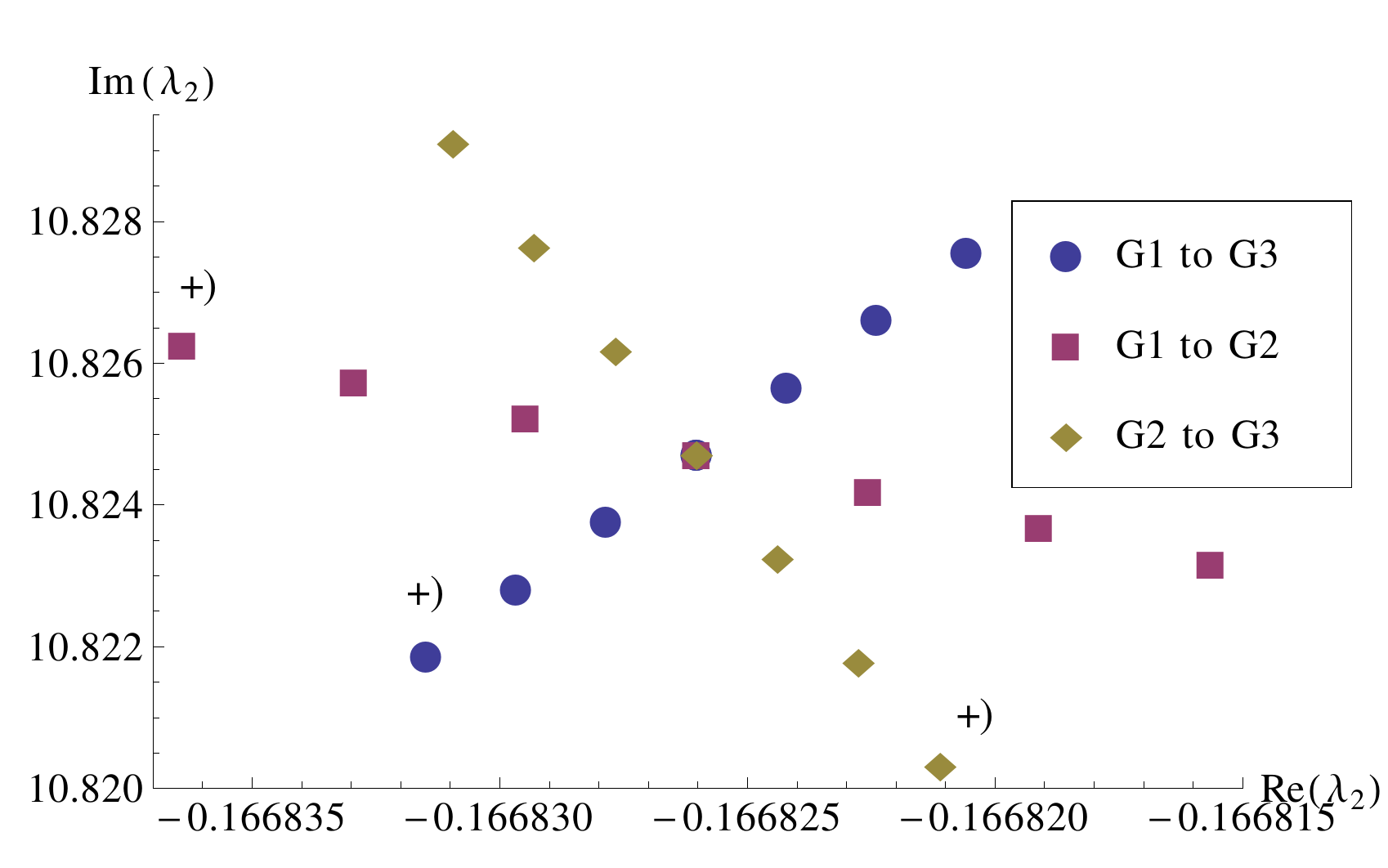}
\caption{\small\label{6BusesEigenvalues} Eigenvalues for redispatches of the 6-bus system.}
\end{figure}
\begin{figure}[h!]
\small
\includegraphics[width=\columnwidth]{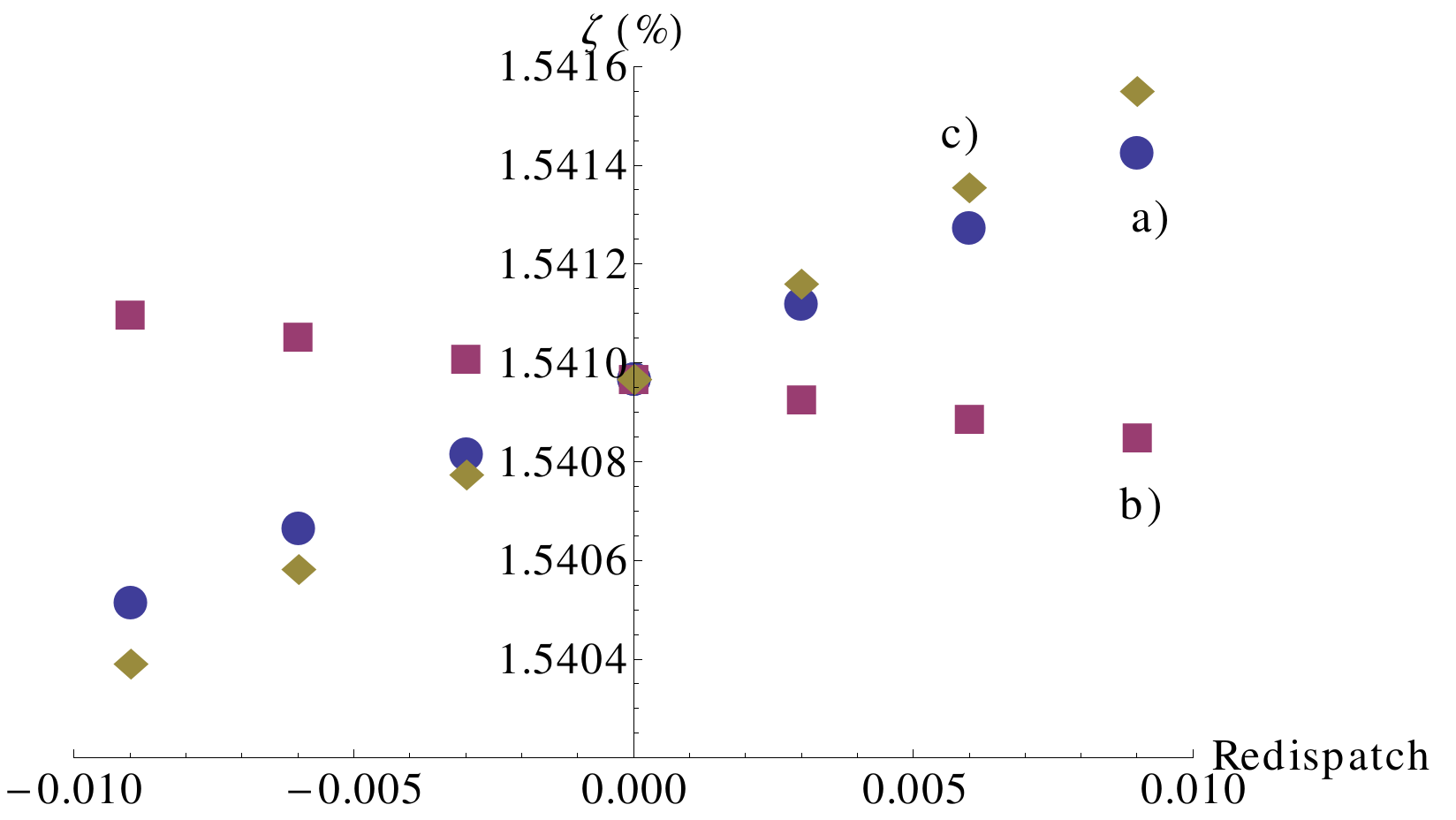}
\caption{\small\label{6BusesDampingPercent}Damping ratio versus redispatch \ \ 
 a) from G1 to G3, 
b) from G1 to G2, c) from G2 to G3.}
\end{figure}

\begin{table}[!h]
\caption{\label{ar and aI 6 Buses Base Case}{Coefficients
$a_{_{rk}}$ and $a_{_{Ik}}$ for the 6-bus system
} }
\centering \begin{tabular}{cccc}
\\
\hline
 $a_{_{r1}}$&- 0.001346  & $a_{_{I1}}$ & - 0.70652\\
$a_{_{r2}}$ & \phantom{-} 0.001275  & $a_{_{I2}}$ & - 1.13594  \\
$a_{_{r3}}$ & \phantom{-} 0.000055  & $a_{_{I3}}$ & - 0.006992  \\
$a_{_{r4}}$ & 0.0 &$a_{_{I4}}$ & - 0.000351\\
$a_{_{r5}}$ & 0.0 & $a_{_{I5}}$& - 0.001029\\
\hline
\end{tabular}
\end{table}
%
From table \ref{ar and aI 6 Buses Base Case}, coefficients related to
the lines 1 and 2 are the biggest components of 
the vectors $a_{_r}$ and $a_{_I}$, but only the coefficients associated
to the line 1 have the same sign, $a_{_{r1}} <0$ and  $a_{_{I_1}}<0$. 
Line 1 connects generator G1, so it is clear from table \ref{ar and aI 6 Buses Base Case}
that increasing G1 helps to damp the oscillation.
Fig.~\ref{6BusesEigenvalues} shows the eigenvalue changes for redispatch
between G1-G3, G1-G2, and G2-G3. When G1 (antinode) increases and G3 (node) 
decreases, $|\sigma_{2}|$ increases and $\omega_{2}$ decreases. If
G1 decreases and G3 increases the effect is opposite. Any other
combination of generators increases or decreases both the real and
imaginary part of $\lambda_2$. 
Table \ref{Critical Eigenvalue and Redispatch Step 6 Buses Base Case}
shows the values of $\lambda_{2f}=d\lambda_2 + \lambda_{2i}$ for different steps of redispatch
between G1 and G3. 
The damping is depicted in Fig.~\ref{6BusesDampingPercent}
as a function of the redispatch of active power. 
The damping ratio improves best when G1 increases and G3 decreases and when
G2 increases and G3 decreases. 

\begin{table}[!h]
\caption{\label{Critical Eigenvalue and Redispatch Step 6 Buses Base Case}{
$\lambda_{2f}$ of redispatch  G1 to G3 in the 6-bus system
} }
\centering \begin{tabular}{cc}\\
Redispatch & $\lambda_{2f}$
\\
\hline
0.009 & -0.166830 + 10.8219j\\ 
0.006 & -0.166830 + 10.8228j\\
0.003 & -0.166828 + 10.8238j\\
0.0 & -0.166826 + 10.8247j \\
-0.003 & -0.166824 + 10.8257j\\
-0.006 & -0.166822 + 10.8266j\\
-0.009 & -0.166821 + 10.8276j\\
\hline
\end{tabular}
\end{table}


\section{Conclusions}

We derive a new formula (\ref{dlambda partial theta partial V}) for the sensitivity of oscillatory eigenvalues 
with respect to generator redispatch.
The motivation is  to understand and improve the damping of interarea oscillations with generator redispatch.

We use a power system dynamic model that expresses both real and reactive power flows and allows for variation of both angle and voltage magnitudes.
The generator dynamics are a simple second order swing equation.
The load modeling allows for frequency dependence and reactive power depending on voltage magnitude, but does not allow 
real power to depend on voltage magnitude.
These modeling assumptions are the usual assumptions permitting energy function analysis of the power system, and 
in particular the network has a symmetric Laplacian. 
Indeed the derivation of the formula exploits the energy function structure.
The hypothesis of the generator dynamic modeling is that there is some equivalent second order  dynamic model for each generator that 
suffices for representing the wide-area oscillations, but that we do not need to know the parameters of each equivalent generator model.
The formula (\ref{dlambda partial theta partial V}) only includes the combined generator dynamics as a common factor that is the same for all redispatches.

In the past, there have multiple unsuccessful attempts to derive a formula with the properties of (\ref{dlambda partial theta partial V}), and sometimes this derivation has been considered to be impossible.
The combination of several ideas in this paper, some new and some old, enables the successful derivation of formula (\ref{dlambda partial theta partial V}):
\begin{enumerate}
\item The new idea of working with the complex symmetric matrix form $x^TQx$ (and not the more obvious  Hermitian matrix form ${\bar x}^TQx$).
\item New ``line" coordinates $(\theta,\nu)$ for the angle differences and logarithm of the product of the voltages across the transmission lines.
These new coordinates greatly simplify parts of the derivation.
\item Quadratic formulation of the eigenvalue problem. 
This formulation was recently applied to a power systems model by Mallada and Tang in \cite{MalladaCDC11}.\footnote{
 \cite{MalladaCDC11} derives the sensitivity of the Fiedler eigenvalue (the smallest magnitude nonzero eigenvalue  of the Laplacian)  
 near saddle node bifurcation  with respect to power injections in the case of constant voltage magnitudes.}
\item The classical assumptions of lossless lines and no dependence of load real power  on voltage magnitude that 
yield the energy function $R$ and a symmetric network Laplacian \cite{BergenHillPAS81,AraposthasisCAS82,NarasimhamurthiCAS84,TsolasCAS85,Pai89,DeMarcoWassnerCCA95}.
\end{enumerate}

The new formula (\ref{dlambda partial theta partial V}) that describes the mode sensitivity 
has a factor $\alpha$ in the denominator that is the same for all 
generator redispatches,
and $\alpha$ depends on the eigenvalue, the equivalent generator dynamics, and the modal eigenvector.
Since the denominator of  (\ref{dlambda partial theta partial V}) is the same for all redispatches, to a large extent we can discriminate the 
effective redispatches by examining the effect of the redispatch on the numerator of  (\ref{dlambda partial theta partial V}).

The numerator of (\ref{dlambda partial theta partial V}) expresses the changes in the mode in terms of 
the changes in angles across lines and load voltage magnitudes caused by the redispatch, with coefficients that
 depend on the mode shape  and the base case power flows in the lines 
and the reactive power load demands.
The base case power flows and the reactive power load demands are available from static state estimation.
The mode shape is  available from synchrophasor measurements, as discussed below.
The new formula (\ref{dlambda partial theta partial V}) is numerically verified in a 10 bus example in section \ref{verify}.

Line coordinates $\theta$ that are the angle differences across the lines are discussed by Bergen and Hill  in  \cite{BergenHillPAS81}.
It is also known that it can be useful to divide the reactive power balance equations by the bus voltage magnitude, and use the 
logarithm  of the bus voltage magnitudes, as, for example, in \cite{OverbyePS94}.  
The  line coordinates $(\theta,\nu)$ are a generalization that includes $\nu$ coordinates 
 that describe the logarithm of the product of the voltage magnitudes  associated with the lines, not the buses.
The line coordinates not only greatly simplify  the derivation of the formula, but are also expected to make the formula easier to interpret when it is applied. There are dependencies between the line coordinates in general meshed networks that are discussed in section \ref{new coordinates}.

The redispatch of real power naturally changes the pattern of real power flows and hence the angles across lines.
Any reactive power flows caused by generator redispatch may also alter the voltage magnitude products across lines.
The numerator of formula (\ref{dlambda partial theta partial V}) identifies in which lines these changes in power flow is most effective.

The main emphasis of this paper is deriving formula (\ref{dlambda partial theta partial V}).
We have also begun to explore the implications and applications of (\ref{dlambda partial theta partial V}) and we now indicate 
some initial conclusions.
\begin{enumerate}
\item In the case that the oscillatory mode has exactly zero
damping, the formula predicts that, to first order, the generator redispatch 
changes only the mode frequency and not the mode damping.
This suggests that generator redispatch could be more effective 
for maintaining sufficient damping than for emergency control 
when damping has vanished.
\item In the special case of considering real power dynamics only with  constant voltage magnitudes,
 the formula (\ref{dlambda partial theta partial V}) reduces to the remarkably simple form (\ref{sensitivity active power case}),
 in which changes in the mode depend on the changes in angles across lines caused by the redispatch, the real power flow in the lines, and the 
line angle coordinates of the mode shape eigenvector $x$.
\item The formula indicates which lines have suitable power flow and eigenvector components to affect oscillation damping.
In particular, it is effective to use the redispatch to change the angle across lines that have both changes in the mode shape across the line 
and sufficient power flow in the right direction.
\end{enumerate}

We note the following considerations and speculations towards implementing  formula (\ref{dlambda partial theta partial V}) to choose the generators to redispatch
that are effective in maintaining suitable oscillation damping or damping ratio.
The complex number $\alpha$ in the denominator of (\ref{dlambda partial theta partial V}) that combines all the equivalent generator dynamics 
is common to all redispatches, so an approximate indication of the argument 
of $\alpha$ is 
probably all that is needed. The base case line power flows are known from the state estimator, and the load flow equations 
can be used to relate the generator redispatches to changes in the angles across lines and the load voltage magnitudes.
The main remaining challenge is to determine the mode shape.

The mode shape is the quadratic eigenvector $x$ corresponding to $\lambda$ and it is easy to obtain  from a conventional right eigenvector.
The mode shape is in principle, and to some considerable extent in practice, available from 
ambient or transient synchrophasor measurements  \cite{TrudnowskiPS08,ChaudhuriPS11,DosiekPS13}.
This is important since it is desirable to use measurements to minimize the use of poorly known dynamic power system models.
Moreover, it is established \cite{PierrePS97,WiesPS03,VanfrettiIREP10,IEEE12} that synchrophasors can make online measurements of the critical eigenvalue  $\lambda$, the oscillatory mode frequency and damping.
And, especially for the low frequency interarea modes, once the mode frequency is known, the mode might have a recurrent and 
fairly robust mode shape.
Then it is conceivable that historical observations or offline computations or general principles about 
the mode shape could be used to augment or interpolate the real-time observations, or that the real time observations 
could be used to verify a predicted mode shape.
Thus some combination of measurements and calculation from models could yield the mode shape
needed to apply the formula to
online  calculations of optimum generation redispatch.

An alternative application of the formula is to 
use it to specify and justify heuristics for oscillation damping based on the mode shape and 
line power flows.  This approach would similarly use a combination of measurements and calculation from models to 
obtain the mode shape, but one might expect that the approximate overall form of the mode shape might suffice.
Our initial results suggest a basis for heuristics for redispatch based on changing the angles across lines with sufficient power flow and 
sufficient changes in the mode shape.
These heuristics would be  similar to heuristics for modal damping due to Fisher and Erlich \cite{FischerIREP,FischerPPT}
that inspired our search for analytic patterns
in modal damping,
and  we would like  to  confirm and refine these heuristics in future work.

More generally, for future work we will 
fully explore the implications and applications of the formula in order to realize its potential for controlling oscillation damping by 
generator redispatch.
The formula could enable some combination of observations, computations
and heuristics to more effectively damp
interarea oscillations.

\section{Acknowledgements}

We gratefully acknowledge support in part from NSF grant CPS-1135825 and 
the Arend J. and Verna V. Sandbulte professorship.  Sarai Mendoza-Armenta gratefully acknowledges
support in part from
Universidad Michoacana de San Nicol\'as de Hidalgo, Conacyt
PhD Scholarship 202024.  Ian Dobson gratefully acknowledge  past 
support towards the solution of this problem 
 coordinated by the Consortium for Electric Reliability Technology Solutions with funding provided in part by the California Energy Commission, Public Interest Energy Research Program, under Work for Others Contract No.~500-99-013, BO-99-2006-P. The Lawrence Berkeley National Laboratory is operated under U.S. Department of Energy Contract No. DE-AC02-05CH11231.
Ian Dobson  thanks Joe Eto for his support of long-term research.

\appendix{\bf Appendix: Jacobian and Quadratic Eigenstructure}
\label{Jacobianquadratic}

In this appendix we show that the eigenvalues and eigenvectors 
of the quadratic form and the Jacobian of the system correspond.
It is convenient to work with the full system of $2n-m$ equations, assumed to have balanced power injections, and  without a reference bus.
Then the system always has a mode with all angles increasing with a zero eigenvalue, which we can neglect.

To compute the eigenvalues of the  
system Jacobian, first we will change 
the  second ordinary 
differential equations (\ref{xnce})  to a set of 
first ordinary differential
equations by defining the variable $\omega_{_{(2n-m)+i}}=\dot{ \delta_i},\ i=1 \ldots m$.
Then the linearized equations become
\begin{align}
\label{linearization Jacobian generators angles}
\dot{\Delta z_i} & = \Delta \omega_{_{(2n-m)+i}}, \ i=1, \ldots m.\\
\label{linearizacion Jacobian angular speed}
\dot{\Delta \omega}_{_{(2n-m)+i}} &= - \frac{d_i}{m_i}\Delta \omega_{_{(2n-m)+i}} -\sum_{j=1}^{2n-m} \frac{L_{ij}}{m_i}\Delta z_j,
\ i=1, \ldots m.\\ 
\label{linearization Jacobian voltages}
 0 & = -\sum_{j=1}^{2n-m} L_{ij}\Delta z_j,  \ i=m+1, \ldots 2n-m.
\end{align}
Writing (\ref{linearization Jacobian generators angles}-\ref{linearization Jacobian voltages}) in matrix form we have 
\begin{align}
\label{Jacobian}
\left(
\begin{matrix}
\dot{\Delta z^{d}} \\
0
\end{matrix}
\right)
=
\left(
\begin{matrix}
J_{11} & J_{12} \\
J_{21} & J_{22} 
\end{matrix}
\right)
\left(
\begin{matrix}
\Delta z^d\\
\Delta z^a\\
\end{matrix}
\right),
\end{align}
where $z^d$ is a vector  of size $2m$ composed by the dynamical variables of the 
system, and $\Delta z^a$ is a vector of size $2(n-m)$ composed of 
the algebraic variables. 

The differential algebraic 
system can be reduced to a purely differential  system by expressing the algebraic 
variables in terms of the dynamic variables  and 
substituting them in the system. This leads to $\Delta z^a=-J^{-1}_{22}J_{21}\Delta z^d$ and the  linearization of the reduction 
\begin{align}
\label{reducedsystem}
\dot{\Delta z^d}=J_{red}{\Delta z^d},
\end{align}
where $J_{red}=J_{11}-J_{12}J^{-1}_{22}J_{21}$. Once
the system is reduced, the symmetry of the Laplacian
of the system is destroyed.

To avoid the reduction (\ref{reducedsystem}), it is better to work directly with 
 the differential-algebraic equations \cite{smed}. (\ref{Jacobian}) 
can be written as a singular ordinary differential equation system
\begin{align}
\label{allsystemc}
E
\left(
\begin{matrix}
\dot{\Delta z^d} \\
\dot{\Delta z^a}
 \end{matrix}
\right)
=
J
\left(
\begin{matrix}
\Delta z^d \\
\Delta z^a
\end{matrix}
\right),
\quad
\text{where}
\quad
E=
\left(
\begin{matrix}
 I & 0 \\
0 & 0
\end{matrix}
\right).
\end{align}
%
To find the eigenvalues associated with 
(\ref{allsystemc}), the generalized
eigenvalue problem has to be solved; i.e., $\mu Ev =\gamma Jv$.
The eigenvalues $\lambda$ are defined as 
$\lambda=\frac{\mu}{\gamma}$. If $\gamma=0$, the eigenvalue $\lambda$
is regarded as infinite.
The infinite eigenvalues arise from  the singularity of the $E$ matrix.

For the finite eigenvalues of
the Jacobian, we can write $Jv= \lambda Ev$.
The eigenvector $v$ is $v=(v^d,v^a)$,
where the size of the vector $v^d$ is the number of dynamics 
variables ($z^d$), and the size
of $v^a$ is the number of algebraic variables. 
It has been proved \cite{smed} that for any triple $(\lambda,v^d,v^a)$
that satisfies (\ref{allsystemc}), the pair $(\lambda,v^d)$
satisfies (\ref{reducedsystem}). Conversely if $(\lambda,v^d)$
satisfies the reduced system, then $(\lambda,v^d,v^a)$ satisfies
the complete system with $v^a=-J^{-1}_{22}J_{21}v^d$, so
the finite eigenvalues of $J$ are the modes
of the system. 

Now we will prove that the finite eigenvalues
of $J$ are finite eigenvalues of $Q$.
Let $v$ be an eigenvector associated with the finite eigenvalue
$\lambda$; that is, $Jv=\lambda E$. Then, from
(\ref{linearization Jacobian generators angles}-\ref{linearization Jacobian voltages}),
\begin{align}
\label{Jacobian eigenvector generators angles0}
\lambda v_i & = v_{(2n-m)+i}\\
\label{linearizacion Jacobian eigenvector angular speed}
\lambda v_{(2n-m)+i} &= - \frac{d_i}{m_i}v_{(2n-m)+i} -\sum_{j=1}^{2n-m} \frac{L_{ij}}{m_i}v_j 
\\ 
\label{linearization Jacobian eigenvector voltages0}
 0 & = -\sum_{j=1}^{2n-m} L_{ij}v_j  .
\end{align}

Using (\ref{Jacobian eigenvector generators angles0}) in (\ref{linearizacion Jacobian eigenvector angular speed}), and multiplying by $m_i$,
\begin{align}
\label{Jacobian eigenvector generators angles}
\lambda^2 m_iv_i + \lambda d_iv_i + \sum_{j=1}^{n+m} L_{ij}v_j &=0 \\ 
\label{linearization Jacobian eigenvector voltages}
\sum_{j=1}^{n+m} L_{ij}v_j  & = 0 .
\end{align}
But (\ref{Jacobian eigenvector generators angles})-(\ref{linearization Jacobian eigenvector voltages})
is (\ref{Qx}). Then the eigenvector $x$ of the quadratic eigenvalue problem
with finite eigenvalue $\lambda$ corresponds exactly to the
eigenvector $v=(\lambda x_g,x)$ of $J$, where $x_g$ is
the vector of components of $x$ corresponding to
the generator angles.


\begin{thebibliography}{99}
\footnotesize
\renewcommand{\em}{\null}



\bibitem{AraposthasisCAS82} A. Arapostathis, S.S. Sastry, P. Varaiya,
IEEE Trans. Circuits and Systems, vol CAS-29, no. 10, October 1982, pp. 673-679.  


\bibitem{BergenHillPAS81}  A.R. Bergen, D.J. Hill, A structure preserving model for power systems stability analysis, IEEE Trans. Power App. Syst., vol. PAS-101, pp. 25-35, Jan. 1981.

\bibitem{ChaudhuriPS11} N.R. Chaudhuri, B. Chaudhuri, Damping and relative mode-shape estimation in near real-time through phasor approach, IEEE Trans. Power Syst., vol. 26, no. 1, pp. 364-373, Feb. 2011.

\bibitem{ChungPS04} C.Y. Chung, L. Wang, F. Howell,  P. Kundur,
Generation rescheduling methods to improve
power transfer capability constrained by
small-signal stability,
IEEE Trans. Power Syst., vol. 19, no. 1, pp. 524-530, Feb. 2004.




\bibitem{Cigre} Cigr\'e Task Force 07 of Advisory Group 01 of Study Committee 38,
Analysis and control of power system oscillations, Paris, December 1996.






\bibitem{DeMarcoWassnerCCA95}     C. L. DeMarco, J. J. Wassner, A generalized
eigenvalue perturbation approach to coherency,  {\em Proc.  IEEE
Conference on
Control Applications}, Albany, NY, September 1995, pp. 605-610.


\bibitem{DiaoPSCE11}
R. Diao, Z. Huang, N. Zhou, Y. Chen, F. Tuffner,  J. Fuller, S. Jin,  J.E Dagle,
Deriving optimal operational rules for mitigating inter-area oscillations,
Power Systems Conference and Exposition, Phoenix AZ USA, March 2011.




\bibitem{DobsonCDC92} I. Dobson, Fernando Alvarado, C. L. DeMarco, Sensitivity
of Hopf bifurcations to power system parameters, Proceedings of the 31st
Conference on Decision and Control, Tucson, Arizona, December 1992.


\bibitem{DobsonPSERC99} I. Dobson, F.L. Alvarado, C.L. DeMarco, P. Sauer, S. Greene, H. Engdahl, J. Zhang,
Avoiding and suppressing oscillations, 
PSerc publication 00-01, December 1999.

\bibitem{DosiekPS13} L. Dosiek, N. Zhou,  J.W. Pierre,  Z. Huang,  D.J. Trudnowski, 
Mode shape estimation algorithms under ambient conditions: A comparative review,
IEEE Transactions on Power Systems, vol. 28, no. 2,
May 2013, pp. 779-787.


\bibitem{FischerIREP} A. Fischer, I. Erlich, Assessment of power
system small signal stability based on mode shape information,
{\em IREP Bulk Power System Dynamics and Control V},
Onomichi, Japan, Aug 2001.

\bibitem{FischerPPT} A. Fischer, I. Erlich, Impact of
long-distance power transits on the dynamic security of
large interconnected power systems,
{\em IEEE Porto Power Tech Conference},
Porto, Portugal, September 2001. 

\bibitem{JonssonPS04} M. Jonsson, M. Begovic, J. Daalder,
A new method suitable for real-time
generator coherency determination,
IEEE Transactions on Power Systems, vol. 19, no. 3,
August 2004, pp. 1473-1482.


\bibitem{mangoPESGM} Z. Huang, N. Zhou, F. Tuffner, Y. Chen, D. Trudnowski, W. Mittelstadt, J. Hauer, J. Dagle,
Improving small signal stability through operating point
adjustment, {\em IEEE PES General Meeting},
Minneapolis, MN USA, July 2010.


\bibitem{mangoReport} Z. Huang, N. Zhou, F.K. Tuffner, Y. Chen, D.J. Trudnowski, MANGO - Modal Analysis for Grid Operation:
A method for damping improvement through operating point adjustment, 
{\em Prepared for the U.S Department of Energy}
October, 2010.





\bibitem{IEEE90} IEEE Power system engineering committee,
{\sl Eigenanalysis and frequency domain 
methods for system dynamic performance}, 
IEEE Publication 90TH0292-3-PWR, 1989.

\bibitem{IEEE94}
IEEE Power Engineering Society Systems Oscillations Working Group,
{\sl Inter-area oscillations in power systems}, IEEE Publication 95~TP~101, 
October 1994.

\bibitem{IEEE12} IEEE Task Force on Identification of Electromechanical Modes,
{\sl Identification of electromechanical modes in power systems},
IEEE Special Publication TP462,
June 2012.


\bibitem{Klein91} M. Klein, G.J. Rogers,  P. Kundur,
A fundamental study of inter-area oscillations in power systems,
{\em IEEE Transactions on Power Systems}, vol. 6,
no. 3, August 1991, pp. 914-921.











\bibitem{MalladaCDC11} E. Mallada, A. Tang, Improving damping of power networks: power scheduling and impedance adaptation,
50th IEEE Conference on Decision and Control and European Control Conference (CDC-ECC),
Orlando, FL, USA, December 2011.

\bibitem{MendozaArmentaPhD} S. Mendoza-Armenta, Analysis of degenerate and interarea oscillations in electric power systems,  (in Spanish),
PhD thesis, Instituto de F\'{\i}sica y Matem\'aticas, Universidad Michoacana, Morelia, Michoac\'an, M\'exico, to appear in 2013.

\bibitem{NamPS00} H.K. Nam, Y.K. Kim, K.S. Shim, K.Y. Lee, A new eigen-sensitivity theory of augmented matrix and its applications to power system stability, IEEE Trans. Power Systems, vol. 15, pp. 363-369, Feb. 2000.

\bibitem{NarasimhamurthiCAS84}
N. Narasimhamurthi and M. T. Musavi, A general energy function for transient stability of power systems, IEEE Trans. Circuits and Systems., vol. CAS-31, pp. 637-645, July 1984.

\bibitem{OverbyePS94}
T.J. Overbye, I. Dobson, C.L. DeMarco, 
Q-V Curve interpretations of energy 
measures for voltage security, 
{\sl IEEE Transactions on Power Systems},
vol. 9, no. 1, Feb. 1994, pp. 331-340.

\bibitem{Pai89}     M. A. Pai, {\em Energy Function Analysis for Power System Stability},
Kluwer Academic Publishers, Boston, 1989.


\bibitem{PierrePS97} J.W. Pierre, D.J. Trudnowski, M.K. Donnelly,
Initial results in electromechanical mode identification from 
ambient data,
{\em IEEE Transactions on Power Systems}, vol. 12,
no. 3, August 1997, pp. 1245-1251.

\bibitem{Rogersbook} G. Rogers, {\sl Power System Oscillations}, Kluwer Academic, 2000.




\bibitem{smed} T. Smed, Feasible eigenvalue sensitivity for large
power systems, IEEE Transactions on Power Systems, Vol. 8, No. 2,
May 1993, pp. 555-563.


\bibitem{TrudnowskiPS08} D.J. Trudnowski, Estimating electromechanical mode shape from synchrophasor measurements, IEEE Trans. Power Syst., vol. 23, no. 3, pp. 1188-1195, Aug. 2008.


\bibitem{TsolasCAS85} N.K. Tsolas, A. Arapostathis, P.P. Varaiya, A structure preserving energy function for power system transient stability analysis,
IEEE Trans.~Circuits and Systems, vol. CAS-32, no. 10, October 1985, pp. 1041-1049.  

\bibitem{VanfrettiIREP10} L. Vanfretti, J.H. Chow,
Analysis of power system oscillations for developing synchrophasor data applications,
2010 IREP Symposium - Bulk Power System Dynamics and Control Ð VIII,  Buzios, Brazil, August 2010,

\bibitem{WangJZU08} Shao-bu Wang, Quan-yuan Jiang, Yi-jia Cao,
WAMS-based monitoring and control of Hopf bifurcations in multi-machine power systems,
Journal of Zhejiang University Science A, vol. 9, no. 6, pp. 840-848, 2008.

\bibitem{WiesPS03} R.W. Wies, J.W. Pierre, D.J. Trudnowski,
Use of ARMA block processing for estimating
stationary low-frequency electromechanical modes
of power systems,
IEEE Trans. Power Syst., vol. 18, no. 1,  Feb. 2003, pp. 167-173.






















































\end{thebibliography}
\end{document}